# STRONG INVARIANCE PRINCIPLES FOR DEPENDENT RANDOM VARIABLES

### By Wei Biao Wu

#### *University of Chicago*

We establish strong invariance principles for sums of stationary and ergodic processes with nearly optimal bounds. Applications to linear and some nonlinear processes are discussed. Strong laws of large numbers and laws of the iterated logarithm are also obtained under easily verifiable conditions.

**1. Introduction.** Strong laws of large numbers (SLLN), laws of the iterated logarithm (LIL), central limit theorems (CLT), strong invariance principles (SIP) and other variants of limit theorems have been extensively studied. Many deep results have been obtained for independent and identically distributed (i.i.d.) random variables. With various weak dependence conditions, some of the results obtained under the i.i.d. assumption have been generalized to dependent random variables. See Ibragimov and Linnik [29], Stout [55], Hall and Heyde [25], Lin and Lu [35], Doukhan [16] and the collection by Eberlein and Taqqu [18] among others.

The primary goal of the paper is to establish a SIP for stationary processes. To this end, we shall develop some moment and maximal inequalities. As our basic tool, a new version of martingale approximation is provided. The martingale method was first applied in Gordin [21] and Gordin and Lifsic [22] and it has undergone substantial improvements. For recent contributions see Merlevède and Peligrad [38], Wu and Woodroofe [67] and Peligrad and Utev [40], where the central limit theory and weak convergence problems are considered. The approximation scheme acts as a bridge which connects stationary processes and martingales. One can then apply results from martingale theory (Chow and Teicher [7] and Hall and Heyde [23]),

---

such as martingale central limit theorems, martingale inequalities, martingale law of the iterated logarithm, martingale embedding theorems, and so on to obtain the desired results for stationary processes.

To implement the martingale method, one needs to know how well stationary processes can be approximated by martingales. In other words, an approximation rate should be obtained. In this paper we shall provide simple sufficient conditions for the existence of $L^q$ $(q > 1)$ martingale approximations as well as approximation rates. A random variable $Z$ is said to be in $L^q$ $(q > 0)$ if $\|Z\|_q := [\mathbb{E}(|Z|^q)]^{1/q} < \infty$. It is convenient to adopt the following formulation. Let $(\xi_i)_{i \in \mathbb{Z}}$ be a stationary and ergodic Markov chain with values in the state space $\mathcal{X}$; let $g : \mathcal{X} \mapsto \mathbb{R}$ be a measurable function for which $\mathbb{E}[g(\xi_0)] = 0$; let the filtration $\mathcal{F}_k = (\ldots, \xi_{k-1}, \xi_k)$, $k \in \mathbb{Z}$. Write $S_n = S_n(g) = \sum_{i=1}^n X_i$, $X_i = g(\xi_i)$. This formulation allows stationary causal processes. Let $\varepsilon_n$ be i.i.d. random elements; let $\xi_n = (\ldots, \varepsilon_{n-1}, \varepsilon_n)$ and $X_n = g(\xi_n)$. Then $(X_n)$ is a causal process and it naturally falls within our framework. As an important category, causal processes have been widely used in practice. Asymptotic results on $S_n$ are useful in the related statistical inference. The $L^q$ martingale approximation, roughly speaking, is to find a martingale $M_n$ with respect to the filter $\mathcal{F}_n$, such that $M_n$ has stationary increments (martingale differences) and the approximation error $\|S_n - M_n\|_q$ is small.

The paper is structured as follows. Section 2 presents an explicit construction of the approximate martingale $M_n$ and a simple and easy-to-use bound of the approximation error $\|S_n - M_n\|_q$. Using those basic tools, we establish in that section various strong laws of large numbers and invariance principles. Some of the results are nearly as sharp as the corresponding ones developed under the i.i.d. assumption. Applications to stationary causal processes are given in Section 3. An interesting feature of the limit theorems in Section 2 is that, besides the necessary moment condition $g(\xi_0) \in L^q$, they basically rely on the magnitude of $\theta_{n,q} = \|\mathbb{E}[g(\xi_n)|\xi_0] - \mathbb{E}[g(\xi_n)|\xi_{-1}]\|_q$, $n \geq 0$. For stationary causal processes $\theta_{n,q}$ is closely related to the *physical and predictive dependence measures* proposed in Wu [64]. It is shown in Section 3 that, for causal processes, those conditions can be easily checked. Applications to linear processes with dependent innovations are discussed in Section 3.2. Proofs are given in Section 4.

We now introduce some notation. Recall $Z \in L^p$ $(p > 0)$ if $\|Z\|_p = [\mathbb{E}(|Z|^p)]^{1/p} < \infty$ and write $\|Z\| = \|Z\|_2$. For ease of reading we list the notation that will be used throughout the paper:

- $(\xi_i)_{i \in \mathbb{Z}}$: a stationary and ergodic Markov chain.
- $S_n = S_n(g) = \sum_{i=1}^n X_i$, where $X_i = g(\xi_i)$ and $g$ is a measurable function.
- $\mathcal{F}_k = (\ldots, \xi_{k-1}, \xi_k)$.
- Projections $\mathcal{P}_k Z = \mathbb{E}(Z|\mathcal{F}_k) - \mathbb{E}(Z|\mathcal{F}_{k-1})$, $Z \in L^1$.

- Let $D_k = \sum_{i=k}^{\infty} \mathcal{P}_k g(\xi_i)$ if the sum converges almost surely.
- $M_k = \sum_{i=1}^{k} D_i$, $R_k = S_k - M_k$, $k \geq 0$.
- $\theta_{n,q} = \|\mathcal{P}_0 g(\xi_n)\|_q$, $n \geq 0$.
- $\Lambda_{n,q} = \sum_{i=0}^{n} \theta_{i,q}$. Let $\theta_{m,q} = 0 = \Lambda_{m,q}$ if $m < 0$.
- Define the tail $\Theta_{m,q} = \sum_{i=m}^{\infty} \theta_{i,q}$ if $\Lambda_{\infty,q} = \lim_{m \to \infty} \Lambda_{m,q} < \infty$.
- $B_q = 18 q^{3/2} (q-1)^{-1/2}$ if $q \in (1,2) \cup (2, \infty)$ and $B_q = 1$ if $q = 2$.
- $\mathbb{B}$: the standard Brownian motion.

By the Markovian property, for $n \geq 0$, $\mathcal{P}_0 g(\xi_n) = \mathbb{E}[g(\xi_n)|\xi_0] - \mathbb{E}[g(\xi_n)|\xi_{-1}]$.

**2. Main results.** Basic tools and some useful moment inequalities are presented in Section 2.1. Sections 2.2 and 2.3 contain LIL, Marcinkiewicz–Zygmund and other forms of laws of large numbers. These results will be useful in proving SIP in Section 2.4.

2.1. *Inequalities and martingale approximations.* Theorem 1 provides moment inequalities and an $L^q$ martingale approximation for $S_n = \sum_{i=1}^{n} g(\xi_i)$. The theorem is proved in Section 4.

THEOREM 1. *Assume that $\mathbb{E}[g(\xi_0)] = 0$ and $g(\xi_0) \in L^q$, $q > 1$. Let $q' = \min(2, q)$. Then* (i)

$$\text{(1)} \qquad \|S_n\|_q^{q'} \leq B_q^{q'} \sum_{i=-n}^{\infty} (\Lambda_{i+n,q} - \Lambda_{i,q})^{q'}.$$

(ii) *Assume additionally that*

$$\text{(2)} \qquad \Theta_{0,q} = \sum_{i=0}^{\infty} \theta_{i,q} < \infty.$$

*Then $D_k := \sum_{i=k}^{\infty} \mathcal{P}_k g(\xi_i)$, $k \in \mathbb{Z}$, are stationary and ergodic $L^q$ martingale differences with respect to $(\mathcal{F}_k)$ and the corresponding martingale $M_k = \sum_{i=1}^{k} D_i$ satisfies*

$$\text{(3)} \qquad \|S_n - M_n\|_q^{q'} \leq 3 B_q^{q'} \sum_{j=1}^{n} \Theta_{j,q}^{q'}.$$

(iii) *Let $S_n^* = \max_{k \leq n} |S_k|$. Then under (2),*

$$\text{(4)} \qquad \|S_n^*\|_q \leq \frac{q B_q}{q-1} n^{1/q'} \Theta_{0,q}.$$

Theorem 1 has two interesting features. First, it provides simple moment bounds for $S_n$. Moment bounds for sums of random variables play an important role in the study of their convergence properties. Second, it presents an explicit construction of approximating martingales as well as the approximation rate (3) under the natural condition (2). The latter condition basically indicates short-range dependence. Hannan [24] proposed condition (2) with $q = 2$ and proved the invariance principle under mixing. Dedecker and Merlevède [11] showed that (2) with $q = 2$ implies the invariance principle without the mixing assumption. McLeish [36] obtained an inequality of type (4) for mixingales with $q = 2$. McLeish's result was improved by Dedecker and Merlevède [10]. With the error term (3), we can quantify the goodness of the approximation and then apply martingale limit theorems. There is a well-developed martingale limit theory and many results established under the independence assumption have their martingale counterparts. These two features are useful in studying the strong convergence of stationary processes.

We now discuss the issue of the uniqueness of the approximation. Observe that (2) and (3) imply $\|S_n - M_n\|_q = o(n^{1/q'})$. It is easily seen that, if $q \geq 2$, then such a martingale construction is unique in the sense that if there is a martingale $M'_n = \sum_{i=1}^n D'_i$ with stationary and ergodic martingale differences $D'_i$ with respect to the filter $\mathcal{F}_i$ such that $\|S_n - M'_n\|_q = o(\sqrt{n})$, then we necessarily have $D_i = D'_i$ almost surely. To this end, note that $\|M'_n - M_n\|_q \leq \|M'_n - S_n\|_q + \|S_n - M_n\|_q = o(\sqrt{n})$. Then $\|M'_n - M_n\| = o(\sqrt{n})$ since $q \geq 2$. Consequently $\|D_1 - D'_1\|\sqrt{n} = o(\sqrt{n})$ and $\|D_1 - D'_1\| = 0$.

REMARK 1. In the case of $q = 2$, the construction of the approximating martingale $(M_k)_{k \geq 0}$ also appeared in Hall and Heyde ([23], page 132), Woodroofe [62] and Volný [60]. To the best of our knowledge, the rate of approximation (3) with a general $q > 1$ is new. Inequality (3) with $q > 2$ is needed in our derivation of the laws of the iterated logarithm and strong invariance principles; see Theorems 2, 3 and 4.

REMARK 2. Yokoyama [71] considered a related problem. The basic assumption imposed in his paper is that there exists a martingale $M_n$ with stationary and ergodic increments such that

$$\|S_n - M_n\|_q = O\{(n \log_2 n)^{q/2}[(\log n) \cdots (\log_m n)^{1+\delta}]^{-1}\}$$

for some $\delta > 0$ and $m \geq 1$. Here $\log_2 n = \log \log n$ and $\log_m n = \log(\log_{m-1} n)$. Under such an assumption, Yokoyama obtained a LIL for $S_n$. For linear processes with i.i.d. innovations, due to the special linearity structure, approximating martingales can be easily constructed; see Section 3 in his paper. However, the latter paper did not address the issue of how to construct approximating martingales for general stationary processes. For nonlinear processes, it is not straightforward to construct such approximating martingales with the desired rate. The explicit construction in Theorem 1(ii) overcomes such restrictions. Some examples are presented in Section 3.

REMARK 3. In the special case $q = 2$, Wu and Woodroofe [68] studied the existence of triangular stationary martingale approximations. Such a martingale approximation scheme means that, for each $n$, there exist stationary martingale differences $D_{n1}, D_{n2}, \ldots$ such that

$$(5) \qquad \max_{1 \le k \le n} \|S_k - M_{nk}\| = o(\|S_n\|),$$

where $M_{nk} = \sum_{i=1}^{k} D_{ni}$. Wu and Woodroofe [68] showed that (5) is equivalent to $\|\mathbb{E}(S_n | \mathcal{F}_0)\| = o(\|S_n\|)$. Under the latter condition, $\|S_n\|$ has the form $\sqrt{nh(n)}$, where $h$ is a slowly varying function. We say that a function $\ell$ is slowly varying if for any $\lambda > 0$, $\lim_{x \to \infty} \ell(\lambda x)/\ell(x) = 1$ (Feller [20], page 275). Their $L^2$ martingale approximation leads to a necessary and sufficient condition for conditional central limit theorems. The martingale approximation (3) differs from (5) in that it allows a general $q > 1$ and that the single-array sequence $(M_k)_{k \ge 0}$ itself is a martingale with stationary increments.

REMARK 4. Inequality (1) is applicable to cases in which $\theta_{n,q}$ decays slowly. For example, if $\theta_{n,q} = (1 + n)^{-\beta}$, $1/q' < \beta < 1$, then (1) gives $\|S_n\|_q = O(n^{1/q'+1-\beta})$. In this case (ii) and (iii) are not applicable since (2) is violated. Corollary 3 contains an application of (1) to strongly dependent processes.

Inequalities concerning $\max_{k \le n} |S_k|$ are important devices in establishing limit theorems for $S_n$. Proposition 1 provides simple maximal inequalities which are useful in the study of strong convergence of stationary processes. Corollary 1 presents an almost sure version of the martingale approximation. Similar maximal inequalities appeared in the literature; see Menchoff [37], Doob [15], Billingsley [5], Serfling [50], Moricz [39] and Lai and Stout [33].

PROPOSITION 1. (i) Let $q > 1$ and $Z_i$, $1 \le i \le 2^d$, be random variables in $\mathcal{L}^q$, where $d$ is a positive integer; let $S_n = Z_1 + \cdots + Z_n$ and $S_n^* = \max_{i \le n} |S_i|$. Then

$$(6) \qquad \|S_{2^d}^*\|_q \le \sum_{r=0}^{d} \left[ \sum_{m=1}^{2^{d-r}} \|S_{2^r m} - S_{2^r(m-1)}\|_q^q \right]^{1/q}.$$

(ii) Let $\{(Y_{k,\phi})_{k \in \mathbb{Z}}, \phi \in \Phi\}$ be a class of centered stationary processes in $L^q, q > 1$; namely for each $\phi \in \Phi$, $(Y_{k,\phi})_{k \in \mathbb{Z}}$ is stationary and $Y_{0,\phi} \in L^q$; let $S_{n,\phi} = Y_{1,\phi} + \cdots + Y_{n,\phi}$. Then

$$(7) \qquad \left\{ \mathbb{E}^* \left[ \max_{k \le n} \sup_{\phi \in \Phi} |S_{k,\phi}|^q \right] \right\}^{1/q} \le \sum_{j=0}^{d} 2^{(d-j)/q} \left\{ \mathbb{E}^* \left[ \sup_{\phi \in \Phi} |S_{2^j,\phi}|^q \right] \right\}^{1/q},$$

where $\mathbb{E}^*$ is the outer expectation: $\mathbb{E}^* Z = \inf\{\mathbb{E} X : X \ge Z, X \text{ is a random variable}\}$.

(iii) *Let $(Z_i)_{i \in \mathbb{Z}}$ be a stationary process and $Z_0 \in L^q$, $q > 0$. Then for all $\delta > 0$,*

$$(8) \qquad \sum_{k=0}^{\infty} \mathbb{P}(S_{2^k}^* \geq 2^{k/q}\delta) \leq 2\delta^{-q}\Delta_q^{q+1}, \qquad where \ \Delta_q = \sum_{j=0}^{\infty}(2^{-j}\|S_{2^j}\|_q^q)^{1/(q+1)}.$$

REMARK 5. Menchoff [37] and Doob [15] considered the special case of (6) with $q = 2$ and uncorrelated random variables. A maximal inequality for stationary processes with $q = 2$ is given in Wu and Woodroofe [68]. The current form (6) allows a general $q > 1$ and nonstationary processes. Inequality (7) is a useful variant of (6) and it is applied in Wu [65].

COROLLARY 1. *Let $g(\xi_1) \in L^q$, $q > 1$, $\mathbb{E}[g(\xi_1)] = 0$ and assume (2). Let $R_k = S_k - M_k$, where $M_k = \sum_{i=1}^{k} D_i$ and $D_k = \sum_{i=k}^{\infty} \mathcal{P}_k g(\xi_i)$. Then $R_n = o_{\text{a.s.}}(n^{1/q})$ under*

$$(9) \qquad \sum_{k=1}^{\infty} k^{-\min[1,(q+4)/(2q+2)]}\Theta_{k,q}^{q/(q+1)} < \infty.$$

PROOF. If $1 < q \leq 2$, since $\Theta_{n,q}$ is nonincreasing in $n$,

$$
\begin{aligned}
\sum_{j=0}^{\infty}\left(2^{-j}\sum_{l=1}^{2^j}\Theta_{l,q}^q\right)^{1/(q+1)} &\leq \sum_{j=0}^{\infty} 2^{-j/(q+1)}\left(\sum_{k=0}^{j} 2^k\Theta_{2^k,q}^q\right)^{1/(q+1)} \\
&\leq \sum_{j=0}^{\infty} 2^{-j/(q+1)}\sum_{k=0}^{j} 2^{k/(q+1)}\Theta_{2^k,q}^{q/(q+1)} \\
&\leq \sum_{k=0}^{\infty}\sum_{j=k}^{\infty} 2^{-j/(q+1)}2^{k/(q+1)}\Theta_{2^k,q}^{q/(q+1)} \\
&= \sum_{k=0}^{\infty} O(\Theta_{2^k,q}^{q/(q+1)}) = \sum_{n=1}^{\infty} n^{-1}O(\Theta_{n,q}^{q/(q+1)}) < \infty.
\end{aligned}
$$

By Theorem 1(ii), $\sum_{j=0}^{\infty}(2^{-j}\|R_{2^j}\|_q^q)^{1/(q+1)} < \infty$. By Proposition 1(iii) and the Borel–Cantelli lemma, we have $R_n = o_{\text{a.s.}}(n^{1/q})$. The other case $q > 2$ similarly follows. □

REMARK 6. If $\Theta_{n,q} = O[n^{\min(0,1/q-1/2)}(\log n)^{\tau}]$, $\tau < -1 - 1/q$, then (9) holds.

REMARK 7. Let $q > 1$. As in Gordin and Lifsic [22], we have $\|R_n\|_q = O(1)$ if

$$(10) \qquad \sum_{i=0}^{\infty} \mathbb{E}[g(\xi_i)|\xi_0] \to H(\xi_0) \qquad (\text{say}) \text{ in } L^q$$

by letting $D_i = H(\xi_i) - \mathbb{E}[H(\xi_i)|\mathcal{F}_{i-1}]$. By Proposition 1(iii) and the Borel–Cantelli lemma, $R_n = o_{\text{a.s.}}(n^{1/q})$. The two sufficient conditions (9) and (10) have different ranges of applicability. Let $g(\xi_n) = \sum_{i=1}^{\infty} i^{-\alpha} \varepsilon_{n-i}$, where $\varepsilon_n \in L^q$ are i.i.d. random variables with $\mathbb{E}(\varepsilon_n) = 0$ and $\mathbb{E}(|\varepsilon_n|^r) = \infty$ for any $r > q$. If $\alpha > 1$, then $\theta_{n,q} = O(n^{-\alpha})$, $\Theta_{n,q} = O(n^{1-\alpha})$ and (9) holds, while (10) requires $\alpha > 1 + 1/q$. On the other hand, however, let $g(\xi_n) = \sum_{i=1}^{\infty} (-1)^i i^{-\beta} \varepsilon_{n-i}$, $1/q < \beta < 1$; then (9) is violated while (10) holds.

2.2. *Laws of large numbers.* The ergodic theorem is probably the best known result on strong convergence of stationary processes. Let $(X_k)$ be a stationary and ergodic process with $\mathbb{E}(|X_0|) < \infty$ and $\bar{X}_n = \sum_{i=1}^{n} X_i/n$; then $\bar{X}_n - \mathbb{E}(X_0) \to 0$ almost surely. However, the convergence of $\bar{X}_n - \mathbb{E}(X_0) \to 0$ can be arbitrarily slow even if $X_k$ is bounded (Krengel [32]). Here we consider the rate of $\bar{X}_n - \mathbb{E}(X_0) \to 0$ under conditions on the decay rates of $\theta_{n,q}$. Corollaries 2 and 3 are easy consequences of Theorem 1 and Proposition 1.

COROLLARY 2. *Let $g(\xi_1) \in L^q$, $q > 1$, $\mathbb{E}[g(\xi_0)] = 0$ and $\ell$ be a positive, nondecreasing slowly varying function.*

(i) *Assume $\Theta_{n,q} = O[(\log n)^{-\alpha}]$, $0 \le \alpha < 1/q$, $q > 2$ and*

$$\sum_{k=1}^{\infty} \frac{k^{-\alpha q}}{[\ell(2^k)]^q} < \infty. \tag{11}$$

*Then $S_n = o_{\text{a.s.}}[\sqrt{n}\ell(n)]$.*

(ii) *Assume (2) with $1 < q \le 2$ and $\sum_{k=1}^{\infty} \ell^{-q}(2^k) < \infty$. Then $S_n = o_{\text{a.s.}}[n^{1/q}\ell(n)]$.*

(iii) *Let $1 < q < 2$ and assume (9). Then $S_n = o_{\text{a.s.}}(n^{1/q})$.*

REMARK 8. Let $\delta > 0$. Then (11) holds if $\ell(n) = (\log n)^{1/q - \alpha} \times (\log n)^{(1+\delta)/q}$. The function $\ell_q(n) = (\log n)^{1/q}(\log \log n)^{(1+\delta)/q}$ satisfies (ii). By (ii), if $\Theta_{0,2} < \infty$, then $S_n = o_{\text{a.s.}}[\sqrt{n}\ell_2(n)]$.

REMARK 9. Let $1 < q < 2$. By the martingale Marcinkiewicz–Zygmund SLLN (cf. Corollary 3 in Tien and Huang [58] or Woyczyński [63]), (10) implies $S_n = o_{\text{a.s.}}(n^{1/q})$; see Remark 7. Dedecker and Merlevède [12] considered Banach-valued random variables.

PROOF OF COROLLARY 2. (i) Let $M_k = \sum_{i=1}^{k} D_i$ and $R_k = S_k - M_k$. By (3) of Theorem 1, $\|R_n\|_q = O[n^{1/2}(\log n)^{-\alpha}]$. By (11) and Proposition 1,

$$\sum_{k=1}^{\infty} \frac{1}{2^{kq/2}\ell^q(2^k)} \mathbb{E}\left[\max_{i \le 2^k} |R_i|^q\right] \le \sum_{k=1}^{\infty} \frac{1}{2^{kq/2}\ell^q(2^k)}\left[\sum_{j=0}^{k} 2^{(k-j)/q}\|R_{2^j}\|_q\right]^q \tag{12}$$

$$= \sum_{k=1}^{\infty} \frac{1}{2^{kq/2}\ell^q(2^k)} O[(2^{k/2}k^{-\alpha})^q] < \infty,$$

which by the Borel–Cantelli lemma implies $\max_{j \le 2^k} |R_j| = o_{a.s.}[2^{k/2}\ell(2^k)]$, and consequently $|R_n| = o_{a.s.}[\sqrt{n}\ell(n)]$ since $\ell$ is slowly varying. By (11) and since $\ell$ is nondecreasing,

$$\sum_{k=1}^{m} \frac{k^{-\alpha q}}{[\ell(2^k)]^q} \ge \sum_{k=1}^{m} \frac{k^{-\alpha q}}{[\ell(2^m)]^q} \sim \frac{m^{1-\alpha q}}{1-\alpha q} \frac{1}{[\ell(2^m)]^q}.$$

So $(\log n)^{1/q-\alpha} = O[\ell(n)]$ and (i) follows from Stout's [54] martingale LIL

$$\limsup_{n \to \infty} \pm \frac{M_n}{\sqrt{2n \log \log n}} = \|D_1\|. \tag{13}$$

(ii) By Theorem 1(iii), $\|S_n^*\|_p^p = O(n)$. So $\sum_{k=1}^{\infty} 2^{-kq/2}\ell^{-q}(2^k)\|S_{2^k}^*\|_p^p < \infty$. By the Borel–Cantelli lemma, since $\ell$ is slowly varying, $S_n = o_{a.s.}[n^{1/q}\ell(n)]$.

(iii) By Corollary 1 and the martingale Marcinkiewicz–Zygmund SLLN, $S_n = R_n + M_n = o_{a.s.}(n^{1/q})$.   $\square$

REMARK 10.   Let $q = 2$ and $\ell^*(n) = \sum_{i=1}^{n} [i\tilde{\ell}(i)]^{-1}$, where $\tilde{\ell}$ is a positive, slowly varying, nondecreasing function. Zhao and Woodroofe [72] proved that $R_n = o_{a.s.}[\sqrt{n\ell^*(n)}]$ under

$$\sum_{i=1}^{\infty} i^{-3/2}\sqrt{\tilde{\ell}(i)}(\log i)\|\mathbb{E}(S_i|\mathcal{F}_0)\| < \infty. \tag{14}$$

Hence $S_n = o_{a.s.}[\sqrt{n\ell^*(n)}]$ if $\log \log n = o[\ell^*(n)]$. Their result does not cover our Corollary 2(ii). On the other hand, however, Corollary 2(ii) does not allow functions like $\ell(n) = \sqrt{\log n}$.

We now give an example where (14) is violated while $\Theta_{0,2} < \infty$. Let $\varepsilon_i$ be i.i.d. with $\mathbb{E}(\varepsilon_i) = 0$ and $\mathbb{E}(\varepsilon_i^2) = 1$ and $X_n = \sum_{i=0}^{\infty} a_i \varepsilon_{n-i}$, where $(a_i)_{i \ge 0}$ are real coefficients. Let $\mathcal{F}_i = (\ldots, \varepsilon_{i-1}, \varepsilon_i)$. Then $\|\mathbb{E}(S_i|\mathcal{F}_0)\|^2 = \sum_{j=0}^{\infty} (a_{j+1} + \cdots + a_{i+j})^2$. Let $a_i = k^{-3/2}$ if $i = 2^k$, $k \in \mathbb{N}$, and $a_i = 0$ if otherwise. Since $a_i$ are nonnegative, $\|\mathbb{E}(S_i|\mathcal{F}_0)\|$ is nondecreasing. Note that $\|\mathbb{E}(S_{2^k}|\mathcal{F}_0)\|^2 \ge 2^k \sum_{l=k}^{\infty} a_{2^l}^2 \ge 2^k k^{-2}/2$. So $\sum_{i=1}^{\infty} i^{-3/2}\|\mathbb{E}(S_i|\mathcal{F}_0)\| = \infty$ since $\sum_{k=1}^{\infty} 2^{-k/2}(2^k k^{-2})^{1/2} = \infty$. Clearly $\Theta_{0,2} < \infty$.

COROLLARY 3.   Let $q > 1$ and $q' = \min(2, q)$.

(i) Assume $\|\mathbb{E}[g(\xi_n)|\mathcal{F}_0]\|_q = O(n^{-\eta})$, $0 < \eta \le 1$. Let $\gamma = \max(1-\eta, 1/q')$.

($i_1$) If $1 - \eta \ne 1/q'$, then $\|S_n^*\|_q = O(n^\gamma)$ and $S_n = o_{a.s.}(n^\gamma(\log n)^{1/q} \times (\log \log n)^{2/q})$.

($i_2$) If $1 - \eta = 1/q'$, then $\|S_n^*\|_q = O(n^\gamma \log n)$ and $S_n = o_{a.s.}(n^\gamma \times (\log n)^{1+1/q}(\log \log n)^{2/q})$.

(ii) Let $\theta_{n,q} = O[n^{-\beta}\tilde{\ell}(n)]$, where $1/q' < \beta < 1$ and $\tilde{\ell}$ is a slowly varying function. Then $\|S_n^*\|_q = O[n^\tau \tilde{\ell}(n)]$ and $S_n = o_{a.s.}[n^\tau \tilde{\ell}(n)\log n]$, $\tau = 1/q' + 1 - \beta$.

Proof. (i) We shall apply the following inequality: there exists a constant $c_q > 0$ such that

$$\|S_{2^d}^*\|_q \le c_q 2^{d/q'} \left[ \|g(\xi_0)\|_q + \sum_{j=0}^d 2^{-j/q'} \|\mathbb{E}(S_{2^j}|\mathcal{F}_0)\|_q \right].$$

Peligrad, Utev and Wu [41] proved the above inequality when $q \ge 2$ [see inequality (10) therein], while Wu and Zhao [70] proved it under $1 < q < 2$. Note that $\|\mathbb{E}(S_{2^j}|\mathcal{F}_0)\|_q \le \sum_{i=1}^{2^j} O(i^{-\eta}) = O(j + 2^{j-j\eta})$.

(i$_1$) Let $1 - \eta \ne 1/q'$. Elementary calculations show that $\|S_{2^d}^*\|_q = O(2^{d\gamma})$. Hence $\|S_n^*\|_q = O(n^\gamma)$. Using the argument in the proof of Corollary 2(ii), the almost sure bound of $S_n$ follows.

(i$_2$) It can be similarly dealt with.

(ii) By Karamata's theorem and Theorem 1(i), $\Lambda_{n,q} = O[n^{1-\beta} \tilde{\ell}(n)]$, and

$$(\|S_n\|_q / B_q)^{q'} \le \sum_{i=-n}^n (\Lambda_{i+n,q} - \Lambda_{i,q})^{q'} + \sum_{i=1+n}^\infty (\Lambda_{i+n,q} - \Lambda_{i,q})^{q'}$$

$$= nO(\Lambda_{2n,q}^{q'}) + \sum_{i=1+n}^\infty O[ni^{1-\beta} \tilde{\ell}(i)]^{q'} = n^{1+(1-\beta)q'} [\tilde{\ell}(n)]^{q'}.$$

By Proposition 1(i), $\|S_{2^d}^*\|_q = O[2^{d\tau} \tilde{\ell}(2^d)]$. By the Borel–Cantelli lemma, (ii) follows. □

REMARK 11. If (11) holds with $\alpha \ge 1/q$, then there exists a law of the iterated logarithm, a more precise form of strong convergence. See Theorem 2(i).

REMARK 12. Let $q > 1$ and $q' = \min(2, q)$. Roughly speaking, since $\Theta_{0,q} < \infty$, Corollary 2 deals with short-range dependent processes. The order of magnitude of $S_n$ is roughly $n^{1/q'}$, up to a multiplicative slowly varying function. Corollary 3 allows long-range dependent processes. For example, if $q = 2$ and $X_n = \sum_{i=0}^\infty a_i \varepsilon_{n-i}$, where $\varepsilon_n$ are i.i.d. with $\mathbb{E}(\varepsilon_n) = 0$, $\mathbb{E}(\varepsilon_n^2) = 1$ and $a_n = n^{-\beta}$, $1/2 < \beta < 1$, then $\|\mathbb{E}(X_n|\mathcal{F}_0)\| = O(n^{1/2-\beta})$ and $S_n = o_{\text{a.s.}}(n^{3/2-\beta} \log n)$.

REMARK 13. With the help of the maximal inequality (6), Corollary 2 can be generalized to nonstationary processes. Let $(X_i)_{i \in \mathbb{N}}$ satisfy $\mathbb{E}(X_i) = 0$ and $\mathbb{E}(|X_i|^q) < \infty$ for all $i \in \mathbb{N}$; let $S_n = X_1 + \cdots + X_n$ for $n \in \mathbb{N}$, $S_0 = 0$ and $S_n^* = \max_{j \le n} |S_j|$. Assume that there exists a nondecreasing function $f$ such that $f(0) = 0$, $\liminf_{n \to \infty} f(n)/f(2n) > 0$ and $\|S_n - S_m\|_q^q \le f(n) - f(m)$ for all integers $n \ge m \ge 0$. By (6), $\|S_{2^d}^*\|_q \le d[f(2^d)]^{1/q}$. If there exists a positive, nondecreasing slowly varying function $\ell(\cdot)$ such that $\sum_{k=1}^{\infty}[k/\ell(2^k)]^q < \infty$, then the argument in the proof of (ii) of Corollary 2 yields $S_n = o_{\text{a.s.}}[f^{1/q}(n)\ell(n)]$. In the special case $q = 2$ and $f(n) = n^\sigma$, $\sigma > 0$, let $\ell(n) = (\log n)^{3/2}(\log\log n)$; then $S_n = o_{\text{a.s.}}[n^{\sigma/2}\ell(n)]$. The latter result is slightly better than the Gaal–Koksma SLLN (Philipp and Stout [44], page 134), which states that $S_n = o_{\text{a.s.}}[n^{\sigma/2}(\log n)^{3/2+\delta}]$ for each $\delta > 0$.

### 2.3. Laws of iterated logarithm.

THEOREM 2. (i) *Assume that* $g(\xi_1) \in L^q$ *for some* $q > 2$, $\mathbb{E}[g(\xi_0)] = 0$ *and*

$$\sum_{i=2}^{\infty}\left[\frac{\Theta_{2^i,q}}{(\log i)^{1/2}}\right]^q < \infty. \tag{15}$$

*Then for* $\sigma := \|\sum_{i=0}^{\infty}\mathcal{P}_0 g(\xi_i)\| < \infty$, *we have, for either choice of sign, that*

$$\limsup_{n\to\infty} \pm\frac{S_n}{\sqrt{2n\log\log n}} = \sigma \tag{16}$$

*almost surely.*

(ii) *Let* $g(\xi_1) \in L^2$ *and* $\mathbb{E}[g(\xi_1)] = 0$. *Then* (16) *holds under* (9) *with* $q = 2$.

LIL for stationary processes has been studied by Philipp [42], Reznik [47], Hall and Heyde [23], Stout [55], Rio [48], Volný and Samek [61] and Zhao and Woodroofe [72] among others. An interesting feature of Theorem 2 is that the conditions (15) and (9) only involve $\theta_{n,q}$ and they are easily verifiable for causal processes; see Proposition 3(ii).

Theorem 2 is proved in Section 4. The key idea is to show that $R_n = o_{\text{a.s.}}(\sqrt{n\log\log n})$. Using the martingale version of Strassen's functional LIL (cf. Heyde and Scott [27] or Basu [3]), we can easily obtain the functional LIL: let $\eta_n(t)$, $t \in [0,1]$, be a function obtained by linearly interpolating $S_i/\sqrt{2n\log\log n}$ at $t = i/n$, $i = 0, 1, \ldots, n$; let $\phi$ be a continuous map from $C[0,1]$ to $\mathbb{R}$. Then under conditions in Theorem 2, $\phi(\eta_n)$ is relatively compact and the set of its limit points coincides with $\phi(F)$, where $F$ is the set of absolutely continuous functions $f$ satisfying $f(0) = 0$ and $\int_0^1 [f'(t)]^2 \le 1$.

Both (i) and (ii) of Theorem 2 concern LIL. They have different ranges of applicability. The former imposes $q > 2$ while in the latter $q = 2$. On the other hand, in some cases, (15) is weaker than (9) with $q = 2$. Let $g(\xi_n) = \sum_{i=0}^{\infty} a_i \varepsilon_{n-i}$, where $\varepsilon_n$ are i.i.d. $L^q$ ($q > 2$) random variables and $a_n = n^{-1}(\log n)^{-\alpha}$, $n \ge 2$, where $\alpha > 1$. Then $\Theta_{n,2} = O[(\log n)^{1-\alpha}]$ and $\Theta_{n,q} = O[(\log n)^{1-\alpha}]$. It is easily seen that (9) needs $\alpha > 5/2$, while (15) only needs $\alpha > 1 + 1/q$.

REMARK 14. Zhao and Woodroofe [72] obtained a LIL by letting $\tilde{\ell}(i) = \log i$ in (14). Their result refines earlier ones by Heyde and Scott [27] and by Heyde [26] in the adapted case. Their LIL and Theorem 2(ii) have different ranges of applicability. Consider the linear process example in Remark 10. Now we let $a_i = k^{-3}$ if $i = 2^k$, $k \in \mathbb{N}$, and $a_i = 0$ if otherwise. So (9) holds with $q = 2$ since $\Theta_{2^k,2} = O(k^{-2})$. Note that $\|\mathbb{E}(S_{2^k}|\mathcal{F}_0)\|^2 \geq 2^k \sum_{l=k}^{\infty} a_{2^l}^2 \geq 2^k k^{-5}/5$. Since $\sum_{k=1}^{\infty} 2^{-k/2} k^{3/2} (2^k k^{-5})^{1/2} = \infty$, (14) is violated [recall $\tilde{\ell}(i) = \log i$]. On the other hand, let $a_i = i^{-1}(\log i)^{-5/2}(\log \log i)^\alpha$, $i \geq 10$. Elementary calculations show that Zhao and Woodroofe's LIL requires $\alpha < -1$ while our Theorem 2(ii) requires $\alpha < -3/2$.

REMARK 15. A simple sufficient condition of (15) is $\Theta_{n,q}^q = O(\log^{-1} n)$, which is weaker than the one in Yokoyama [71] even in the special case of linear processes. Let $g(\xi_n) = \sum_{i=0}^{\infty} a_i \varepsilon_{n-i}$, where $\varepsilon_n$ are i.i.d. and $a_n$ are real coefficients. For $a_n = O[n^{-1}(\log n)^{-\alpha}]$ and $\varepsilon_0 \in L^q$, $q > 2$, Yokoyama's result requires $\alpha \geq 1 + 2/q$ whereas $\alpha \geq 1 + 1/q$ is enough to ensure (15).

2.4. *Strong invariance principles.* The strong invariance principle studied here means the almost sure approximation of partial sums of random variables by Brownian motions. With such approximations, asymptotic properties of partial sums can be obtained from those of Brownian motions. Strong invariance principles are quite useful in statistical inference of time series and have received considerable attention in probability theory. Strassen [56, 57] initiated the study for i.i.d. random variables and stationary and ergodic martingale differences. Komlós et al. [30, 31] considered approximating sums of i.i.d. random variables by Brownian motions and obtained optimal results. Motivated by the usefulness of such approximation scheme, researchers have established many results concerning mixingales and strongly mixing processes of various types; see Berkes and Philipp [4], Bradley [6], Eberlein [17], Shao [51], Rio [48], Lin and Lu [35], Dedecker and Prieur [13] and Philipp and Stout [44] among others. Philipp [43] gave an excellent review.

Here we shall apply the moment and maximal inequalities and SLLN developed in Sections 2.1 and 2.2 to study strong approximations of partial sums of stationary processes. To state such results, one often needs to enlarge the underlying probability space and redefine the stationary process without changing its distribution. For brevity we simply say that there is a richer probability space and a standard Brownian motion $\mathbb{B}$ such that the partial sum process can be approximated by $\mathbb{B}$ with certain rates.

It seems that our method is quite effective and it leads to nearly optimal approximation rates. Let $\chi_q(n) = n^{1/q}(\log n)^{1/2}$ if $2 < q < 4$ and $\chi_q(n) = n^{1/4}(\log n)^{1/2}(\log \log n)^{1/4}$ if $q \geq 4$; let $\iota_q(n) = n^{1/q}(\log n)^{1/2+1/q}(\log \log n)^{2/q}$. Recall Theorem 1 for $D_k = \sum_{i=k}^{\infty} \mathcal{P}_k g(\xi_i)$.

THEOREM 3.   *Let $g(\xi_0) \in L^q$, $q > 2$, $\mathbb{E}[g(\xi_0)] = 0$ and $\sigma = \|D_k\|$. Let $q^* = \min(q, 4)$.*

(i)  *Assume $\Theta_{n,q^*} = O[n^{1/q^*-1/2}(\log n)^{-1}]$  and*

$$(17) \qquad\qquad \sum_{k=1}^{\infty} \|\mathbb{E}(D_k^2|\mathcal{F}_0) - \sigma^2\|_{q^*/2} < \infty.$$

*Then on a richer probability space, there exists a standard Brownian motion $\mathbb{B}$ such that*

$$(18) \qquad\qquad |S_n - \mathbb{B}(\sigma^2 n)| = O_{\mathrm{a.s.}}[\chi_q(n)].$$

(ii)  *Assume $\Theta_{n,q^*} = O(n^{1/q^*-1/2})$  and*

$$(19) \qquad\qquad \sum_{k=1}^{\infty} \|\mathcal{P}_0(D_k^2)\|_{q^*/2} < \infty.$$

*Then on a richer probability space, there exists a standard Brownian motion $\mathbb{B}$ such that*

$$(20) \qquad\qquad |S_n - \mathbb{B}(\sigma^2 n)| = O_{\mathrm{a.s.}}[\iota_{q^*}(n)].$$

REMARK 16.   In Theorem 3, conditions (17) and (19) involve the decays of $\|\mathbb{E}(D_k^2|\mathcal{F}_0) - \sigma^2\|_{q^*/2}$ and $\|\mathcal{P}_0(D_k^2)\|_{q^*/2}$. Interestingly, it turns out that for stationary causal processes both quantities have simple and easy-to-use bounds; see Proposition 3 in Section 3.

THEOREM 4.   *Let $g(\xi_0) \in L^q$, $2 < q \le 4$, and $\mathbb{E}[g(\xi_0)] = 0$. Assume $\Theta_{n,q} = O(n^{-\eta})$  and $\|\mathbb{E}(D_n^2|\mathcal{F}_0) - \sigma^2\|_{q/2} = O(n^{-\eta})$, $\eta > 0$. Then on a richer probability space, there exists a standard Brownian motion $\mathbb{B}$ such that $|S_n - \mathbb{B}(\sigma^2 n)| = o_{\mathrm{a.s.}}[n^{\gamma/2}(\log n)^{3/2}]$, where $\gamma = \max(1 - \eta, 2/q)$.*

Theorems 3 and 4 give explicit approximation rates. It is interesting to note that, when $2 < q \le 4$, the rates (18) and (20) are optimal up to multiplicative logarithmic factors. Komlós et al. [30, 31] showed that, if $(X_k)_{k \in \mathbb{Z}}$ are i.i.d. random variables with mean 0 and $q$th moment, $q > 2$, then the approximation rate is $o(n^{1/q})$ and it cannot be improved to $o(n^{\delta})$ for any $\delta < 1/q$. Strong invariance principles with the rates $n^{1/2-\delta}$ or $o[(n \log \log n)^{1/2}]$ are widely considered in the literature; see Philipp and Stout [44], Philipp [43] and Eberlein [17].

SIP plays an important role in statistical inference. For change-point and trend analysis see Csörgő and Horvath [8]. Wu and Zhao [69] considered statistical inference of trends in time series and constructed simultaneous confidence bands for mean trends with asymptotically correct nominal coverage probabilities. Shao and Wu [52] studied asymptotic properties of the local Whittle estimator of the Hurst (long-memory) parameter for fractionally integrated nonlinear time series models. In the latter two papers we applied the SIP of the form (20) and the rate $O_{\mathrm{a.s.}}[\iota_{q^*}(n)]$ is needed to control certain error terms.

It is unclear whether a bound of the form $O[n^{1/q}\ell(n)]$, where $\ell$ is a slowly varying function, can be obtained when $q > 4$. The arguments in Komlós et al. [30, 31], which heavily depend on the independence assumption, cannot be directly applied to stationary processes. In our proof of Theorem 3 we apply Strassen's [57] martingale embedding method. The best bound that Strassen's method can result in is $\chi_4(n)$. In the special case of linear processes, since the approximate martingales are sums of i.i.d. random variables, sharp bounds can be derived.

PROPOSITION 2. *Let $X_n = \sum_{i=0}^{\infty} a_i \varepsilon_{n-i}$, where $\varepsilon_n \in L^q$, $q > 2$, are i.i.d. with mean 0 and the real coefficients $a_n$ are square summable. Assume that $A_j := \sum_{i=j}^{\infty} a_i$, $j \geq 0$, exists and*

$$(21) \qquad \sum_{j=1}^{\infty} (2^{-j}\Xi_{2j}^q)^{1/(q+1)} < \infty \qquad \text{where } \Xi_n = \sqrt{\sum_{i=1}^{n} A_i^2 + \sum_{i=n+1}^{\infty} (A_i - A_{i-n})^2}.$$

*Then on a richer probability space, there exists a standard Brownian motion $\mathbb{B}$ such that*

$$(22) \qquad |S_n - \mathbb{B}(\sigma^2 n)| = o_{\mathrm{a.s.}}(n^{1/q}) \qquad \text{where } \sigma = |A_0| \|\varepsilon_0\|.$$

*In particular, (21) holds under (9), or $\sum_{i=1}^{\infty} A_i^2 < \infty$, or*

$$(23) \qquad \sum_{i=1}^{\infty} \sqrt{\sum_{j=i}^{\infty} a_j^2} < \infty.$$

PROOF. Let $M_k = A_0 \sum_{i=1}^{k} \varepsilon_i$. Then $-R_n = M_n - S_n = \sum_{i=1}^{\infty} (A_i - A_{i-n} \times \mathbf{1}_{i>n}) \varepsilon_{n+1-i}$. By Rosenthal's inequality, $\|R_n\|_q \leq C \|R_n\|$ for some constant $C < \infty$. Since $\|R_n\| = \Xi_n \|\varepsilon_0\|$, by Corollary 1, (21) implies $R_n = o_{\mathrm{a.s.}}(n^{1/q})$ and hence (22) by the Hungarian construction of Komlós et al. Clearly $\sum_{i=1}^{\infty} A_i^2 < \infty$ implies $\Xi_n = O(1)$. Under (23), $\sum_{n=0}^{\infty} \mathbb{E}(X_n|\mathcal{F}_0)$ converges in $L^2$, we have as in Remark 7 that $\|R_n\| = O(1)$ and then (21). $\square$

Phillips [45] studied the consistency problem of log-periodogram regression. A key tool in the latter paper is a SIP for linear processes with the rate $o_{\mathrm{a.s.}}(n^\zeta)$, $\zeta \in (2/q, 1/2)$ under

$$(24) \qquad \sum_{i=1}^{\infty} i|a_i| < \infty$$

and $\varepsilon_i \in L^q$, $q > 4$. We now compare Phillips' result with our SIP (22). First, our result requires $q > 2$ while $q > 4$ is needed in Phillips' result. Second, Phillips' condition (24) implies $\Theta_{n,q} = O(n^{-1})$ and hence (9). Our conditions are much weaker. Third, our bound $o_{\text{a.s.}}(n^{1/q})$ is optimal and is sharper than Phillips' $o_{\text{a.s.}}(n^\zeta)$. Akonom [1] obtained a bound $o_{\mathbb{P}}(n^{1/q})$ under (24).

**3. Applications.** Let $(\varepsilon_n)_{n \in \mathbb{Z}}$ be i.i.d. random variables and let $X_n = g(\xi_n)$, where $\xi_n = (\ldots, \varepsilon_{n-1}, \varepsilon_n)$ and $g$ is a measurable function such that $X_n$ is a proper random variable. We can view the random variables $\varepsilon_i, i \leq n$, as the input to a system, $g$ as a filter and $X_n = g(\xi_n)$ as the output of the system. The class of causal processes is huge and it includes a variety of nonlinear time series models (Priestley [46], Tong [59] and Stine [53]).

In this section we shall apply the results in Section 2 to stationary causal processes. Due to the structure of $\xi_n$, the filter $\mathcal{F}_n$ can be naturally defined as the sigma algebra generated by $\xi_n$. Write $\mathbb{E}(Z|\mathcal{F}_n) = \mathbb{E}(Z|\xi_n)$. To apply the moment inequalities and limit theorems presented in Section 2, one needs to verify the conditions therein [cf. (2), (9), (15), (17) and (19)] which are based on the quantities $\|\mathcal{P}_0 g(\xi_n)\|_q$, $\|\mathbb{E}[g(\xi_n)|\mathcal{F}_0]\|_q$, $\|\mathcal{P}_0(D_n^2)\|_q$ and $\|\mathbb{E}(D_n^2|\mathcal{F}_0) - \sigma^2\|_q$. It turns out that for causal processes such quantities can be effectively handled by the coupling method. Let $(\varepsilon_n')_{n \in \mathbb{Z}}$ be an i.i.d. copy of $(\varepsilon_n)_{n \in \mathbb{Z}}$ and let $\xi_k' = (\ldots, \varepsilon_{n-1}', \varepsilon_n')$. For $k \geq 0$ let $\xi_k^* = (\xi_0', \varepsilon_1, \ldots, \varepsilon_{k-1}, \varepsilon_k)$ and $\tilde{\xi}_k = (\xi_{-1}, \varepsilon_0', \varepsilon_1, \ldots, \varepsilon_{k-1}, \varepsilon_k)$. In $\xi_k^*$ the whole "past" $\xi_0$ is replaced by the i.i.d. copy $\xi_0'$, while $\tilde{\xi}_k$ only couples a single innovation $\varepsilon_0$. Let $h_m(\xi_k) = \mathbb{E}[g(\xi_{k+m})|\xi_k]$ be the $m$-step conditional expectation, $m \geq 1$, and write $h(\xi_k) = h_1(\xi_k)$. For $k \geq 0$ define

$$
\begin{aligned}
(25) \qquad &\tilde{\alpha}_k = \|h(\xi_k) - h(\tilde{\xi}_k)\|_q, \qquad \alpha_k^* = \|h(\xi_k) - h(\xi_k^*)\|_q, \\
&\tilde{\beta}_k = \|g(\xi_k) - g(\tilde{\xi}_k)\|_q, \qquad \beta_k^* = \|g(\xi_k) - g(\xi_k^*)\|_q.
\end{aligned}
$$

The input/output viewpoint provides another way to look at the dependence. Wu [64] introduced the *physical dependence measure* $\tilde{\beta}_k$, which measures how much the process will deviate, measured by the $L^q$ distance, from the original orbit $(g(\xi_k))_{k \geq 0}$ if we change the current input $\varepsilon_0$ to an i.i.d. copy $\varepsilon_0'$. By (27), $\|\mathcal{P}_0 g(\xi_n)\|_q$ has the same order of magnitude as $\|h_n(\xi_0) - h_n(\tilde{\xi}_0)\|_q$. Note that $h_n(\xi_0) = \mathbb{E}[g(\xi_n)|\xi_0]$ is the $n$-step ahead predicted mean. So $\omega_n = \|h_n(\xi_0) - h_n(\tilde{\xi}_0)\|_q$ measures the contribution of $\varepsilon_0$ in predicting future values. In Wu [64] $\omega_n$ is called *predictive dependence measure*. The short-range dependence condition $\Theta_{0,q} < \infty$ means that the cumulative contribution of $\varepsilon_0$ in predicting future values is finite. Dedecker and Doukhan [9] proposed $\mathbb{E}|g(\xi_k) - g(\ldots, 0, \varepsilon_1, \ldots, \varepsilon_k)|$, a similar version of $\beta_k^*$. See Wu [64] for more discussions on input/output dependence measures.

Proposition 3 provides bounds for $\|\mathcal{P}_0 g(\xi_n)\|_q$, $\|\mathbb{E}[g(\xi_n)|\mathcal{F}_0]\|_q$, $\|\mathcal{P}_0(D_n^2)\|_q$ and $\|\mathbb{E}(D_n^2|\mathcal{F}_0) - \sigma^2\|_q$ based on the quantities in (25). Despite the fact that

the martingale differences $D_n$ are constructed from $g(\xi_i)$ in a complicated manner, there exist simple bounds for $\|\mathcal{P}_0(D_n^2)\|_q$ and $\|\mathbb{E}(D_n^2|\mathcal{F}_0) - \sigma^2\|_q$.

PROPOSITION 3. *Let* $g(\xi_0) \in L^q$, $q > 1$, $q' = \min(2, q)$, *and* $\mathbb{E}[g(\xi_0)] = 0$.

(i) *Let* $k \geq 0$. *Then* $\tilde{\alpha}_k \leq 2\tilde{\beta}_{k+1}$ *and* $\alpha_k^* \leq 2\beta_{k+1}^*$.

(ii) *For* $n \geq 1$ *we have*

$$(26) \qquad \tfrac{1}{2}\|h_n(\xi_0) - h_n(\xi_0^*)\|_q \leq \|\mathbb{E}(g(\xi_n)|\xi_0)\|_q$$

$$\leq \min[\|h_n(\xi_0) - h_n(\xi_0^*)\|_q, \alpha_{n-1}^*]$$

*and*

$$(27) \qquad \tfrac{1}{2}\|h_n(\xi_0) - h_n(\tilde{\xi}_0)\|_q \leq \|\mathcal{P}_0 g(\xi_n)\|_q \leq \min[\|h_n(\xi_0) - h_n(\tilde{\xi}_0)\|_q, \tilde{\alpha}_{n-1}].$$

(iii) *Let* $q > 2$ *and* $\sum_{i=1}^{\infty} \tilde{\alpha}_i < \infty$. *Then* $D_k = \sum_{i=k}^{\infty} \mathcal{P}_k g(\xi_i) \in L^q$. *Additionally, for* $k \in \mathbb{N}$,

$$(28) \qquad \|\mathbb{E}(D_k^2|\xi_0) - \sigma^2\|_{q/2} \leq 8c_q\beta_k^* + 8c_q \sum_{i=k}^{\infty} \min(\alpha_i^*, \tilde{\alpha}_{i-k}),$$

*where* $\sigma = \|D_k\|$ *and* $c_q = \|D_k\|_q$, *and*

$$(29) \qquad \|\mathcal{P}_0(D_k^2)\|_{q/2} \leq 8c_q\tilde{\beta}_k + 8c_q \sum_{i=k}^{\infty} \tilde{\alpha}_i.$$

REMARK 17. In Proposition 3(iii), the results are expressed in terms of the one-step conditional expectation $h(\xi_k) = \mathbb{E}(g(\xi_{k+1})|\xi_k)$. It is straightforward to generalize them to the $m$-step conditional expectations $h_m(\xi_k) = \mathbb{E}[g(\xi_{k+m})|\xi_k]$.

COROLLARY 4. *Let* $g(\xi_0) \in L^q$, $2 < q \leq 4$, *and* $\mathbb{E}[g(\xi_0)] = 0$. *Assume*

$$(30) \qquad \sum_{k=1}^{\infty} (\tilde{\beta}_k + k\tilde{\alpha}_k) < \infty.$$

*Then* (20) *holds. In particular,* (30) *holds if*

$$(31) \qquad \sum_{k=1}^{\infty} k\tilde{\beta}_k < \infty.$$

Corollary 4 is an easy consequence of (ii) of Theorem 3 and Proposition 3 in view of $\Theta_{n,q} \leq \sum_{k=n}^{\infty} \tilde{\alpha}_{k-1} = O(n^{-1})$, by (27) and (30). By Proposition 3(i), $\tilde{\alpha}_k \leq 2\tilde{\beta}_{k+1}$, so (31) implies (30). Corollary 4 provides simple sufficient conditions for strong invariance principles. The quantities $\tilde{\alpha}_k$ and $\tilde{\beta}_k$ are directly related to the data-generating mechanisms.

3.1. *Transforms of linear processes.* Let $X_n = \sum_{i=0}^{\infty} a_i \varepsilon_{n-i}$, where $\varepsilon_k$ are i.i.d., $\varepsilon_k \in L^q$, $q > 1$, and $(a_n)_{n \geq 0}$ is a real sequence. Let $K$ be a Lipschitz continuous function; namely, there exists a constant $C < \infty$ such that $|K(x) - K(y)| \leq C|x - y|$ for all $x, y \in \mathbb{R}$. Let $g(\xi_n) = K(X_n) - \mathbb{E}[K(X_n)]$. Then $\theta_{n,q} \leq \beta_n = O(\|a_n(\varepsilon_0 - \varepsilon_0')\|_q) = O(|a_n|)$. So (31) holds if $\sum_{i=1}^{\infty} i|a_i| < \infty$.

3.2. *Linear processes with dependent innovations.* Let $(\varepsilon_k)_{k \in \mathbb{Z}}$ be i.i.d. random elements and let $(a_n)_{n \geq 1}$ be a sequence of real numbers; let

$$(32) \qquad X_n = \sum_{i=0}^{\infty} a_i \eta_{n-i} \qquad \text{where } \eta_n = G(\ldots, \varepsilon_{n-1}, \varepsilon_n)$$

and $G$ is a measurable function. Linear processes with dependent innovations have attracted considerable attention recently. Models of such type have been proposed in econometric time series analysis. See Baillie, Chung and Tieslau [2], Romano and Thombs [49], Hauser and Kunst [25], Lien and Tse [34] and Wu and Min [66] among others. Appropriate conditions on the innovations are needed for the asymptotic analysis of $(X_n)$. Let $q > 2$. Assume that there exist $C > 0$ and $r \in (0, 1)$ such that for all $n \in \mathbb{N}$,

$$(33) \qquad \mathbb{E}[|G(\xi_0) - G(\xi_n^*)|^q] \leq Cr^n,$$

where $\xi_n^* = (\xi_0', \varepsilon_1, \ldots, \varepsilon_n)$. The property (33) is called the geometric-moment contraction condition (Hsing and Wu [28], Wu and Shao [67]) and it is satisfied for many nonlinear time series. In particular, for the iterated random function (Elton [19] and Diaconis and Freedman [14]) $\eta_n = R(\eta_{n-1}, \varepsilon_n)$, let

$$L_\varepsilon = \sup_{x \neq x'} \frac{|R(x, \varepsilon) - R(x', \varepsilon)|}{|x - x'|}$$

be the Lipschitz coefficient. Then (33) is satisfied if $\mathbb{E}(L_\varepsilon^q) < 1$ and $\mathbb{E}[|R(x_0, \varepsilon)|^q] < \infty$ for some $x_0$. Wu and Min [66] showed that (33) holds for GARCH processes.

COROLLARY 5. *Let $2 < q \leq 4$ and assume (33).*

(i) *Assume $\sum_{i=n}^{\infty} |a_i| = O(\log^{-1/q} n)$. Then the LIL (16) holds.*
(ii) *Assume $\sum_{i=n}^{\infty} |a_i| = O(n^{1/q-1/2}/\log n)$. Then the SIP (18) holds.*

**4. Proofs.**

PROOF OF THEOREM 1. (i) Consider two cases $1 < q \leq 2$ and $q > 2$ separately. For $1 < q \leq 2$, since $\mathcal{P}_k S_n$, $k = -\infty, \ldots, n-1, n$, form martingale differences, by Burkholder's inequality,

$$(\|S_n\|_q/B_q)^q \leq \mathbb{E}\left(\sum_{k=-\infty}^n |\mathcal{P}_k S_n|^2\right)^{q/2} \leq \sum_{k=-\infty}^n \mathbb{E}(|\mathcal{P}_k S_n|^q)$$

$$\leq \sum_{k=-\infty}^n (\Lambda_{n-k,q} - \Lambda_{-k,q})^q,$$

where we have applied the inequalities $(|a_1| + |a_2| + \cdots)^{q/2} \leq |a_1|^{q/2} + |a_2|^{q/2} + \cdots$ and $\|\mathcal{P}_k S_n\|_q \leq \Lambda_{n-k,q} - \Lambda_{-k,q}$. Clearly $B_q$ can be chosen to be 1 if $q = 2$. If $q > 2$, then $\|\cdot\|_{q/2}$ is a norm and by Burkholder's inequality

$$(\|S_n\|_q/B_q)^2 \leq \left\|\sum_{k=-\infty}^n |\mathcal{P}_k S_n|^2\right\|_{q/2} \leq \sum_{k=-\infty}^n \|(\mathcal{P}_k S_n)^2\|_{q/2}$$

$$= \sum_{k=-\infty}^n \|\mathcal{P}_k S_n\|_q^2 \leq \sum_{k=-\infty}^n (\Lambda_{n-k,q} - \Lambda_{-k,q})^2.$$

(ii) By (2), $D_k \in L^q$. Let $R_k = S_k - M_k$. By the triangle inequality, for $1 \leq j \leq n$, $\|\mathcal{P}_j R_n\|_q \leq \Theta_{n+1-j,q}$ and, for $j \leq 0$, $\|\mathcal{P}_j R_n\|_q \leq \Theta_{1-j,q} - \Theta_{n+1-j,q}$. We first consider the case $q \geq 2$. Since $R_n = \sum_{j=-\infty}^n \mathcal{P}_j R_n$ and the summands form martingale differences and $\mathcal{P}_j \mathcal{P}_k = 0$ if $j \neq k$, again by Burkholder's inequality,

$$(34) \qquad \mathbb{E}(|R_n|^q) \leq B_q^q \mathbb{E}\left[\sum_{j=-\infty}^n (\mathcal{P}_j R_n)^2\right]^{q/2},$$

which together with the triangle inequality $\|\sum_{j=-\infty}^n Z_j\|_{q/2} \leq \sum_{j=-\infty}^n \|Z_j\|_{q/2}$ yields

$$\frac{\|R_n\|_q^2}{B_q^2} \leq \sum_{j=-\infty}^0 \|\mathcal{P}_j R_n\|_q^2 + \sum_{j=1}^n \|\mathcal{P}_j R_n\|_q^2$$

$$\leq \sum_{j=-\infty}^{-n} (\Theta_{1-j,q} - \Theta_{n+1-j,q})^2$$

$$\quad + \sum_{j=-n+1}^0 (\Theta_{1-j,q} - \Theta_{n+1-j,q})^2 + \sum_{j=1}^n \Theta_{n+1-j,q}^2$$

$$\leq \sum_{j=-\infty}^{-n} (\Theta_{1-j,q} - \Theta_{n+1-j,q})\Theta_n + 2\sum_{j=1}^n \Theta_{n+1-j,q}^2$$

$$\leq n\Theta_{n,q}^2 + 2\sum_{j=1}^n \Theta_{n+1-j,q}^2 \leq 3\sum_{j=1}^n \Theta_{j,q}^2.$$

The other case in which $1 < q < 2$ is similar and it follows from (34) and

$$\mathbb{E}\left[\sum_{j=-\infty}^{n} (\mathcal{P}_j R_n)^2\right]^{q/2} \leq \sum_{j=-\infty}^{n} \mathbb{E}[(\mathcal{P}_j R_n)^2]^{q/2} \leq \sum_{j=-\infty}^{n} \mathbb{E}(|\mathcal{P}_j R_n|^q)$$

$$\leq \sum_{j=-\infty}^{0} (\Theta_{1-j,q} - \Theta_{n+1-j,q})^q + \sum_{j=1}^{n} \Theta_{n+1-j,q}^q$$

$$\leq 3\sum_{j=1}^{n} \Theta_{j,q}^q.$$

Combining the two cases, we have (3).

(iii) For $j \geq 0$ let $Z_{j,n} = \sum_{k=1}^{n} \mathcal{P}_{k-j} g(\xi_k)$. Let $p = q/(q-1)$. We apply McLeish's [36] argument. Again we consider two cases $1 < q \leq 2$ and $q > 2$ separately. Let $q > 2$. Notice that for a fixed $j$, $\mathcal{P}_{k-j} g(\xi_k)$ form stationary martingale differences in $k = 1, \ldots, n$. By Burkholder's inequality and the triangle inequality,

$$(35) \quad \begin{aligned} \|Z_{j,n}\|_q^2 &\leq B_q^2 \left\| \sum_{k=1}^{n} [\mathcal{P}_{k-j} g(\xi_k)]^2 \right\|_{q/2} \\ &\leq B_q^2 \sum_{k=1}^{n} \|[\mathcal{P}_{k-j} g(\xi_k)]^2\|_{q/2} = B_q^2 n \|\mathcal{P}_0 g(\xi_j)\|_q^2. \end{aligned}$$

If $1 < q \leq 2$, then

$$(36) \quad \mathbb{E}\left[\sum_{k=1}^{n} |\mathcal{P}_{k-j} g(\xi_k)|^2\right]^{q/2} \leq \sum_{k=1}^{n} \mathbb{E}[|\mathcal{P}_{k-j} g(\xi_k)|^2]^{q/2} = n \|\mathcal{P}_0 g(\xi_k)\|_q^q.$$

By Doob's inequality (cf. Hall and Heyde [23], Theorem 2.2, page 15), $\|\max_{k \leq n} |Z_{j,k}|\|_q \leq p\|Z_{j,n}\|_q$. By (35) and (36), $\|Z_{j,n}\|_q \leq B_q n^{1/q'} \|\mathcal{P}_0 g(\xi_j)\|_q$. Since $\tilde{S}_k = \sum_{j=0}^{\infty} Z_{j,k}$ and

$$\|S_n^*\|_q \leq \sum_{j=0}^{\infty} \left\| \max_{k \leq n} |Z_{j,k}| \right\|_q \leq \sum_{j=0}^{\infty} p B_q n^{1/q'} \|\mathcal{P}_0 g(\xi_j)\|_q = p B_q n^{1/q'} \Theta_{0,q},$$

we have (4). $\quad \square$

PROOF OF PROPOSITION 1. (i) Let $p = q/(q-1)$ and $\Lambda = \sum_{r=0}^{d} \lambda_r^{-p}$, where

$$\lambda_r = \left[\sum_{m=1}^{2^{d-r}} \|S_{2^r m} - S_{2^r(m-1)}\|_q^q\right]^{-1/(p+q)}.$$

For the positive integer $k \leq 2^d$, write its dyadic expansion $k = 2^{r_1} + \cdots + 2^{r_j}$, where $0 \leq r_j < \cdots < r_1 \leq d$, and $k(i) = 2^{r_1} + \cdots + 2^{r_i}$. By Hölder's inequality,

$$
\begin{aligned}
|S_k|^q &\leq \left[ \sum_{i=1}^{j} |S_{k(i)} - S_{k(i-1)}| \right]^q \\
&\leq \left[ \sum_{i=1}^{j} \lambda_{r_i}^{-p} \right]^{q/p} \left[ \sum_{i=1}^{j} \lambda_{r_i}^{q} |S_{k(i)} - S_{k(i-1)}|^q \right] \\
&\leq \Lambda^{q/p} \sum_{i=1}^{j} \lambda_{r_i}^{q} \sum_{m=1}^{2^{d-r_i}} |S_{2^{r_i}m} - S_{2^{r_i}(m-1)}|^q \\
&\leq \Lambda^{q/p} \sum_{r=0}^{d} \lambda_{r}^{q} \sum_{m=1}^{2^{d-r}} |S_{2^r m} - S_{2^r(m-1)}|^q,
\end{aligned}
$$

which entails $\|S_{2^d}^*\|_q^q \leq \Lambda^{q/p} \sum_{r=0}^{d} \lambda_r^q \lambda_r^{-p-q} = \Lambda^q$ and hence (6).

(ii) It easily follows from the argument in (i) since $\sup_\theta |S_{k,\theta}| \leq \sum_{i=1}^{j} \sup_\theta |S_{k(i),\theta} - S_{k(i-1),\theta}|$.

(iii) It suffices to prove (8) with $\Delta_q < \infty$. Let $r_j = (2^{-j} \|S_{2^j}\|_q^q)^{1/(q+1)}/\Delta_q$, $j \geq 0$. Then $\sum_{j=0}^{\infty} r_j = 1$. By the argument in (i), we have

$$
\begin{aligned}
\sum_{d=0}^{\infty} \mathbb{P}(S_{2^d}^* \geq 2^{d/q}\delta) &\leq \sum_{d=0}^{\infty} \sum_{j=0}^{d} \mathbb{P}\left[ \max_{m \leq 2^{d-j}} |S_{2^j m} - S_{2^j(m-1)}| \geq 2^{k/q}\delta r_j \right] \\
&\leq \sum_{j=0}^{\infty} \sum_{d=j}^{\infty} 2^{d-j} \mathbb{P}(|S_{2^j}| \geq 2^{d/q}\delta r_j) \\
&\leq \sum_{j=0}^{\infty} \sum_{k=0}^{\infty} 2^k \mathbb{P}[2^{-j}|S_{2^j}|^q (\delta r_j)^{-q} \geq 2^k] \\
&\leq \sum_{j=0}^{\infty} 2\mathbb{E}[2^{-j}|S_{2^j}|^q (\delta r_j)^{-q}] = 2\delta^{-q}\Delta_q^{q+1},
\end{aligned}
$$

where the inequality $\sum_{k=0}^{\infty} 2^k \mathbb{P}(|Z| \geq 2^k) \leq 2\mathbb{E}(|Z|)$ is applied. $\quad\square$

PROOF OF THEOREM 2. (i) Recall $D_k = \sum_{i=k}^{\infty} \mathcal{P}_k g(\mathcal{F}_i)$, $M_k = \sum_{i=1}^{k} D_i$ and $R_k = S_k - M_k$. Let $p = q/(q-1)$, $\alpha = 1/2 - 1/q > 0$, $\lambda_j = 2^{\alpha j/2}$ and $\psi_k = \Theta_{2^k, q}$. Note that $\Theta_{n,q}$ is nonincreasing in $n$. By (3) of Theorem 1,

$$
\|R_{2^j}\|_q \leq C\left( \sum_{m=1}^{2^j} \Theta_m^2 \right)^{1/2} \leq C\left( \sum_{i=0}^{j} 2^i \psi_i^2 \right)^{1/2} \leq C \sum_{i=0}^{j} 2^{i/2} \psi_i,
$$

where $C$ is a constant that does not depend on $j$ and it may vary from line to line. So

$$\sum_{j=0}^{k} 2^{(k-j)/q} \|R_{2^j}\|_q \leq C \sum_{j=0}^{k} 2^{(k-j)/q} \sum_{i=0}^{j} 2^{i/2} \psi_i$$

$$\leq C \sum_{i=0}^{k} 2^{i/2} \psi_i \sum_{j=i}^{\infty} 2^{(k-j)/q}$$

$$\leq C 2^{k/q} \sum_{i=0}^{k} 2^{\alpha i} \psi_i.$$

By Hölder's inequality,

$$\sum_{i=0}^{k} 2^{\alpha i} \psi_i \leq \left( \sum_{i=0}^{k} \lambda_i^q \psi_i^q \right)^{1/q} \left( \sum_{i=0}^{k} 2^{\alpha i p} \lambda_i^{-p} \right)^{1/p} \leq \left( \sum_{i=0}^{k} \lambda_i^q \psi_i^q \right)^{1/q} C \lambda_k.$$

Hence by (15) and Proposition 1,

$$\sum_{k=3}^{\infty} \frac{1}{\sqrt{2^k \log k^q}} \mathbb{E}\left[ \max_{i \leq 2^k} |R_i|^q \right] \leq \sum_{k=3}^{\infty} \frac{1}{\sqrt{2^k \log k^q}} \left[ \sum_{j=0}^{k} 2^{(k-j)/q} \|R_{2^j}\|_q \right]^q$$

$$\leq C \sum_{k=3}^{\infty} \frac{2^{-\alpha k q}}{(\log k)^{q/2}} \left( \sum_{i=0}^{k} \lambda_i^q \psi_i^q \right) C \lambda_k^q$$

$$\leq C \sum_{i=0}^{\infty} \lambda_i^q \psi_i^q \sum_{k=\max(i,3)}^{\infty} \frac{2^{-\alpha k q}}{(\log k)^{q/2}}$$

$$\leq C + C \sum_{i=3}^{\infty} \lambda_i^q \psi_i^q \frac{2^{-\alpha i q}}{(\log i)^{q/2}} < \infty.$$

By the Borel–Cantelli lemma, $\max_{i \leq 2^k} |R_i| = o_{\text{a.s.}}(\sqrt{2^k \log k})$, which completes the proof by Stout's martingale LIL.

(ii) By Corollary 1, $R_n = o_{\text{a.s.}}(\sqrt{n})$. So the LIL easily follows. □

PROOF OF THEOREM 3. (i) As in Theorem 1, let $M_k = \sum_{i=1}^{k} D_i$ and $R_k = S_k - M_k$. If $2 < q < 4$, by (3) of Theorem 1, $\|R_n\| = O[n^{1/q}(\log n)^{-1}]$. By Proposition 1,

$$\left\| \max_{k \leq 2^d} |R_k| \right\|_q \leq \sum_{j=0}^{d} 2^{(d-j)/q} \|R_{2^j}\|_q$$

$$\leq \sum_{j=0}^{d} 2^{(d-j)/q} O\left[ 1 + \sum_{i=2}^{2^j} i^{2/q-1} (\log i)^{-2} \right]^{1/2} = O(2^{d/q} \log d).$$

Hence we have $\max_{k \leq 2^d} |R_k| = o_{\text{a.s.}}[\chi_q(2^d)]$ and, by the argument following (12), $R_n = o_{\text{a.s.}}[\chi_q(n)]$. The case in which $q \geq 4$ similarly follows.

Now we shall deal with the martingale part and show that (18) holds with $S_n$ replaced by $M_n$. To this end, we apply Strassen's [57] martingale embedding method. Let

$$J_n = \sum_{k=1}^{n} [\mathbb{E}(D_k^2|\mathcal{F}_{k-1}) - \sigma^2].$$

If $2 < q < 4$, let $g(\mathcal{F}_k) = \mathbb{E}(D_{k+1}^2|\mathcal{F}_k) - \sigma^2$, then $g(\mathcal{F}_k) \in L^{q/2}$ and (17) implies $J_n = o_{\text{a.s.}}(n^{2/q})$; see Remark 7. By the martingale version of the Skorokhod representation theorem (Strassen [57], Hall and Heyde [23]), on a possibly richer probability space, there exist a standard Brownian motion $\mathbb{B}$ and nonnegative random variables $\tau_1, \tau_2, \ldots$ with partial sums $T_k = \sum_{i=1}^{k} \tau_i$ such that, for $k \geq 1$, $M_k = \mathbb{B}(T_k)$, $\mathbb{E}(\tau_k|\mathcal{F}_{k-1}) = \mathbb{E}(D_k^2|\mathcal{F}_{k-1})$ and

$$(37) \qquad \mathbb{E}(\tau_k^{q/2}|\mathcal{F}_{k-1}) \leq C_q \mathbb{E}(|D_k|^q|\mathcal{F}_{k-1})$$

almost surely, where $C_q$ is a constant which only depends on $q$. Let $Q_n = \sum_{k=1}^{n} [\tau_k - \mathbb{E}(\tau_k|\mathcal{F}_{k-1})]$. Since $\mathbb{E}(\tau^{q/2}) < \infty$ and $1 < q/2 < 2$, by the martingale Marcinkiewicz–Zygmund SLLN, $Q_n = o_{\text{a.s.}}(n^{2/q})$. Notice that

$$T_n - n\sigma^2 = Q_n + \sum_{k=1}^{n} [\mathbb{E}(\tau_k|\mathcal{F}_{k-1}) - \mathbb{E}(D_k^2|\mathcal{F}_{k-1})] + J_n.$$

Therefore $T_n - n\sigma^2 = o_{\text{a.s.}}(n^{2/q})$, and

$$(38) \qquad \begin{aligned} \max_{k \leq n} |\mathbb{B}(T_k) - \mathbb{B}(\sigma^2 k)| &\leq \max_{k \leq n} \sup_{x\,:\,|x - \sigma^2 k| \leq n^{2/q}} |\mathbb{B}(x) - \mathbb{B}(\sigma^2 k)| \\ &= o_{\text{a.s.}}[n^{1/q}(\log n)^{1/2}]. \end{aligned}$$

If $q \geq 4$, by (ii) of Theorem 2, (17) implies the LIL

$$\limsup_{n \to \infty} \pm \frac{J_n}{t_n} = \left\| \sum_{k=0}^{\infty} \mathcal{P}_0 \mathbb{E}(D_k^2|\mathcal{F}_{k-1}) \right\| = \left\| \sum_{k=1}^{\infty} \mathcal{P}_0 D_k^2 \right\| < \infty,$$

where $t_n = \sqrt{2n \log \log n}$. By the martingale LIL (Stout [54]), $Q_n = O_{\text{a.s.}}(t_n)$. So there is a constant $C < \infty$ such that $|T_n - n\sigma^2| \leq Ct_n$ almost surely and similarly

$$\begin{aligned} \max_{k \leq n} |\mathbb{B}(T_k) - \mathbb{B}(\sigma^2 k)| &\leq \max_{k \leq n} \sup_{x\,:\,|x - \sigma^2 k| \leq Ct_n} |\mathbb{B}(x) - \mathbb{B}(\sigma^2 k)| \\ &= O[t_n^{1/2}(\log n)^{1/2}] = O[\chi_4(n)] \end{aligned}$$

holds almost surely.

(ii) As in (12), we have by Proposition 1 that the condition $\Theta_{n,q^*} = O(n^{1/q^*-1/2})$ implies $R_n = o_{a.s.}[\iota_{q^*}(n)]$. By Corollary 2(ii), let $\ell(n) = (\log n)^{2/q^*}(\log\log n)^{(2+\delta)/q^*}$ (see also Remark 8), (19) entails $J_n = o_{a.s.}[n^{2/q^*} \times \ell(n)]$, which yields the desired result in view of the proof of (i).  $\square$

PROOF OF THEOREM 4.    Without loss of generality let $0 < \eta \le 1 - 2/q$. By (3), $\|R_n\|_q^2 = (n^{1-2\eta} + \log n)$, which entails $R_n = o_{a.s.}(n^{\gamma/2})$ by the argument in (12). We now deal with $M_n$. By Corollary 3(i), the condition $\|\mathbb{E}(D_n^2|\mathcal{F}_0) - \sigma^2\|_{q/2} = O(n^{-\eta})$ implies

$$J_n = \sum_{k=1}^n [\mathbb{E}(D_k^2|\mathcal{F}_{k-1}) - \sigma^2] = o_{a.s.}[n^\gamma(\log n)^{1+2/q}(\log\log n)^{4/q}].$$

So the desired result follows from the argument of (38) in the proof of Theorem 3.  $\square$

PROOF OF PROPOSITION 3.    (i) Let $i \ge k+1$. Then we have the identity

$$
\begin{aligned}
(39)\qquad \mathbb{E}[h(\tilde{\xi}_{i-1})|\xi_k] &= \mathbb{E}[h(\tilde{\xi}_{i-1})|\xi_{-1}, \varepsilon_j, 1 \le j \le k] \\
&= \mathbb{E}[h(\xi_{i-1})|\xi_{-1}, \varepsilon_j, 1 \le j \le k] = \mathbb{E}[h(\xi_{i-1})|\tilde{\xi}_k].
\end{aligned}
$$

Hence by Jensen's inequality,

$$
\begin{aligned}
\tilde{\alpha}_k &= \|\mathbb{E}[g(\xi_{k+1})|\xi_k] - \mathbb{E}[g(\tilde{\xi}_{k+1})|\tilde{\xi}_k]\|_q \\
&\le \|\mathbb{E}[g(\xi_{k+1})|\xi_k] - \mathbb{E}[g(\tilde{\xi}_{k+1})|\xi_k]\|_q + \|\mathbb{E}[g(\tilde{\xi}_{k+1})|\xi_k] - \mathbb{E}[g(\tilde{\xi}_{k+1})|\tilde{\xi}_k]\|_q \\
&\le \|g(\xi_{k+1}) - g(\tilde{\xi}_{k+1})\|_q + \|\mathbb{E}[g(\tilde{\xi}_{k+1})|\tilde{\xi}_k] - \mathbb{E}[g(\tilde{\xi}_{k+1})|\tilde{\xi}_k]\|_q \le 2\tilde{\beta}_{k+1}.
\end{aligned}
$$

The inequality $\alpha_k^* \le 2\beta_{k+1}^*$ can be similarly proved.

(ii) We only consider (27) since (26) can be similarly established. Observe that $\mathbb{E}[g(\xi_n)|\xi_0] = \mathbb{E}[h(\xi_{n-1})|\xi_0]$ and $\mathbb{E}[g(\xi_n)|\xi_{-1}] = \mathbb{E}[h(\xi_{n-1})|\xi_0]$. By Jensen's inequality, $\|\mathcal{P}_0 g(\xi_n)\|_q = \|\mathbb{E}[h(\xi_{n-1}) - h(\tilde{\xi}_{n-1})|\xi_0]\|_q \le \tilde{\alpha}_{n-1}$. Since $\mathbb{E}[h_n(\tilde{\xi}_0)|\xi_0] = \mathbb{E}[h_n(\xi_0)|\xi_{-1}]$,

$$\|\mathcal{P}_0 g(\xi_n)\|_q = \|h_n(\xi_0) - \mathbb{E}[h_n(\tilde{\xi}_0)|\xi_0]\|_q \le \|h_n(\xi_0) - h_n(\tilde{\xi}_0)\|_q.$$

On the other hand,

$$
\begin{aligned}
\|h_n(\xi_0) - h_n(\tilde{\xi}_0)\|_q &\le \|h_n(\xi_0) - h_{n+1}(\xi_{-1})\|_q + \|h_{n+1}(\xi_{-1}) - h_n(\tilde{\xi}_0)\|_q \\
&= 2\|\mathcal{P}_0 g(\xi_n)\|_q.
\end{aligned}
$$

Combining the preceding inequalities, we have (27).

(iii) Under the proposed conditions, $D_k = \sum_{i=k}^{\infty} \mathcal{P}_k g(\xi_i) \in L^q$ by (i). Write $D_k = D(\xi_k)$ and $\mathcal{P}_k g(\xi_i) = H_{i-k}(\xi_k)$, where $D(\cdot)$ and $H_{i-k}(\cdot)$ are measurable functions. We first show (29). Note that for $k \geq 1$, $\mathbb{E}(D_k^2|\xi_{-1}) = \mathbb{E}[D(\tilde{\xi}_k)^2|\xi_0]$. Since $q/2 > 1$, by Jensen's and Schwarz's inequalities,

$$\|\mathcal{P}_0(D_k^2)\|_{q/2} = \|\mathbb{E}[D(\xi_k)^2 - D(\tilde{\xi}_k)^2|\xi_0]\|_{q/2}$$
$$\leq \|D(\xi_k)^2 - D(\tilde{\xi}_k)^2\|_{q/2}$$
$$\leq \|D(\xi_k) - D(\tilde{\xi}_k)\|_q \|D(\xi_k) + D(\tilde{\xi}_k)\|_q$$
$$\leq 2c_q\|D(\xi_k) - D(\tilde{\xi}_k)\|_q.$$

Notice that $D(\xi_k) = \sum_{i=k}^{\infty} H_{i-k}(\xi_k)$. Then

$$\|D(\xi_k) - D(\tilde{\xi}_k)\|_q \leq \|H_0(\xi_k) - H_0(\tilde{\xi}_k)\|_q$$
$$+ \sum_{i=k+1}^{\infty} \|H_{i-k}(\xi_k) - H_{i-k}(\tilde{\xi}_k)\|_q. \tag{40}$$

Clearly, we have by (i) that

$$\|H_0(\xi_k) - H_0(\tilde{\xi}_k)\|_q \leq \|g(\xi_k) - g(\tilde{\xi}_k)\|_q + \|h(\xi_{k-1}) - h(\tilde{\xi}_{k-1})\|_q$$
$$= \tilde{\beta}_k + \tilde{\alpha}_{k-1} \leq 3\tilde{\beta}_k.$$

To show (29), it remains to verify that $\|H_{i-k}(\xi_k) - H_{i-k}(\tilde{\xi}_k)\|_q \leq 4\tilde{\alpha}_{i-1}$ holds for $i \geq k+1$. To this end, since $H_{i-k}(\xi_k) = \mathbb{E}[h(\xi_{i-1})|\xi_k] - \mathbb{E}[h(\xi_{i-1})|\xi_{k-1}]$ and $H_{i-k}(\tilde{\xi}_k) = \mathbb{E}[h(\tilde{\xi}_{i-1})|\tilde{\xi}_k] - \mathbb{E}[h(\tilde{\xi}_{i-1})|\tilde{\xi}_{k-1}]$,

$$\|H_{i-k}(\xi_k) - H_{i-k}(\tilde{\xi}_k)\|_q \leq \|\mathbb{E}[h(\xi_{i-1})|\xi_k] - \mathbb{E}[h(\tilde{\xi}_{i-1})|\tilde{\xi}_k]\|_q$$
$$+ \|\mathbb{E}[h(\xi_{i-1})|\xi_{k-1}] - \mathbb{E}[h(\tilde{\xi}_{i-1})|\tilde{\xi}_{k-1}]\|_q. \tag{41}$$

Since $i \geq k+1$, by (39) and Jensen's inequality, we have $\mathbb{E}[h(\tilde{\xi}_{i-1})|\xi_k] = \mathbb{E}[h(\xi_{i-1})|\tilde{\xi}_k]$ and

$$\|\mathbb{E}[h(\xi_{i-1})|\xi_k] - \mathbb{E}[h(\tilde{\xi}_{i-1})|\tilde{\xi}_k]\|_q$$
$$\leq \|\mathbb{E}[h(\xi_{i-1})|\xi_k] - \mathbb{E}[h(\tilde{\xi}_{i-1})|\xi_k]\|_q$$
$$+ \|\mathbb{E}[h(\tilde{\xi}_{i-1})|\xi_k] - \mathbb{E}[h(\tilde{\xi}_{i-1})|\tilde{\xi}_k]\|_q$$
$$\leq \tilde{\alpha}_{i-1} + \|\mathbb{E}[h(\xi_{i-1})|\tilde{\xi}_k] - \mathbb{E}[h(\tilde{\xi}_{i-1})|\tilde{\xi}_k]\|_q$$
$$\leq 2\tilde{\alpha}_{i-1}. \tag{42}$$

Similarly, $\|\mathbb{E}[h(\xi_{i-1})|\xi_{k-1}] - \mathbb{E}[h(\tilde{\xi}_{i-1})|\tilde{\xi}_{k-1}]\| \leq 3\tilde{\alpha}_{i-1}$. By (41), $\|H_{i-k}(\xi_k) - H_{i-k}(\tilde{\xi}_k)\|_q \leq 4\tilde{\alpha}_{i-1}$ and (29) follows.

The same argument also applies to (28). We need to replace $\tilde{\xi}_i$, $\tilde{\alpha}_i$ and $\tilde{\beta}_i$ by the $*$ versions $\xi_i^*$, $\alpha_i^*$ and $\beta_k^*$. Then $\|H_0(\xi_k) - H_0(\xi_k^*)\|_q \leq 2\beta_k^*$ and $\|H_{i-k}(\xi_k) - H_{i-k}(\xi_k^*)\|_q \leq 4\alpha_{i-1}^*$ hold for $i \geq k+1$. Notice that we also have

$$\|H_{i-k}(\xi_k) - H_{i-k}(\xi_k^*)\|_q \leq \|H_{i-k}(\xi_k)\|_q + \|H_{i-k}(\xi_k^*)\|_q \leq 2\tilde{\alpha}_{i-k-1}.$$

So (28) follows from (40).    □

PROOF OF COROLLARY 5.    (i) By Theorem 2, it suffices to verify (15). Observe that (33) implies $\|\mathcal{P}_0 \eta_n\|_q \leq Cr^n$. Therefore,

$$
\begin{aligned}
\sum_{m=n}^{\infty} \|\mathcal{P}_0 X_m\|_q &\leq \sum_{m=n}^{\infty} \sum_{i=0}^{m} |a_{m-i}| Cr^i \leq C \sum_{i=0}^{\infty} \sum_{m=\max(i,n)}^{\infty} |a_{m-i}| r^i \\
&= C \sum_{i=0}^{n-1} \sum_{m=n}^{\infty} |a_{m-i}| r^i + C \sum_{i=n}^{\infty} \sum_{m=i}^{\infty} |a_{m-i}| r^i \\
&= \sum_{i=0}^{n-1} r^i O[\log^{-1/q}(n-i)] + \sum_{i=n}^{\infty} O(r^i) = O(\log^{-1/q} n),
\end{aligned}
$$

(43)

which implies $\Theta_{n,q}^q = O(\log^{-1} n)$ and consequently (15).

(ii) Under the condition on $(a_n)$, we have $\Theta_{n,q^*} = O[n^{1/q^*-1/2}(\log n)^{-1}]$ in view of the argument in (43). So the corollary follows from (i) of Theorem 3 and Lemma 1.    □

LEMMA 1.    *For the process* (32), *assume* $\sum_{i=0}^{\infty} |a_i| < \infty$ *and* (33), $2 < q \leq 4$. *Then* (17) *holds.*

PROOF.    For $n \geq k$, let $\mathcal{P}_k \eta_n = L_{n-k}(\xi_k)$. Similarly as in (42), we have

$$
\begin{aligned}
\|L_{n-k}(\xi_k) &- L_{n-k}(\xi_k^*)\|_q \\
&\leq \|\mathbb{E}(\eta_n|\xi_k) - \mathbb{E}(\eta_n^*|\xi_k^*)\|_q + \|\mathbb{E}(\eta_n|\xi_{k-1}) - \mathbb{E}(\eta_n^*|\xi_{k-1}^*)\|_q \\
&\leq 4\|\eta_n - \eta_n^*\|_q = O(r^{n/q})
\end{aligned}
$$

in view of (33). Let $\rho = r^{1/q}$. Since $\mathcal{P}_k g(\xi_n) = \sum_{j=0}^{\infty} a_j \mathcal{P}_k \eta_{n-j} = \sum_{j=0}^{n-k} a_j \times L_{n-j-k}(\xi_k)$ and $D_k = D(\xi_k) = \sum_{n=k}^{\infty} \mathcal{P}_k g(\xi_n)$,

$$
\begin{aligned}
\sum_{k=1}^{\infty} \|D(\xi_k) - D(\xi_k^*)\|_q &\leq \sum_{k=1}^{\infty} \sum_{n=k}^{\infty} \left\| \sum_{j=0}^{n-k} a_j [L_{n-j-k}(\xi_k) - L_{n-j-k}(\xi_k^*)] \right\|_q \\
&\leq \sum_{k=1}^{\infty} \sum_{n=k}^{\infty} \sum_{j=0}^{n-k} |a_j| \rho^{n-j} \\
&\leq \sum_{k=1}^{\infty} \rho^k (1-\rho)^{-1} \sum_{j=0}^{\infty} |a_j| < \infty,
\end{aligned}
$$

which entails (17) since $\|\mathbb{E}[D^2(\xi_k)|\xi_0] - \sigma^2\|_{q/2} \leq 2\|D_k\|_q \|D(\xi_k) - D(\xi_k^*)\|_q$. $\square$

**Acknowledgments.** I am grateful to the two referees and an Associate Editor for their many helpful comments. I would like to thank Jan Mielniczuk for valuable suggestions.

## REFERENCES

[1] AKONOM, J. (1993). Comportement asymptotique du temps d'occupation du processus des sommes partielles. *Ann. Inst. H. Poincaré Probab. Statist.* **29** 57–81. MR1204518

[2] BAILLIE, R. T., CHUNG, C. F. and TIESLAU, M. A. (1996). Analysing inflation by the fractionally integrated ARFIMA–GARCH model. *J. Appl. Econometrics* **11** 23–40.

[3] BASU, A. K. (1973). A note on Strassen's version of the law of the iterated logarithm. *Proc. Amer. Math. Soc.* **41** 596–601. MR0329007

[4] BERKES, I. and PHILIPP, W. (1979). Approximation theorems for independent and weakly dependent random vectors. *Ann. Probab.* **7** 29–54. MR0515811

[5] BILLINGSLEY, P. (1968). *Convergence of Probability Measures.* Wiley, New York. MR0233396

[6] BRADLEY, R. C. (1983). Approximation theorems for strongly mixing random variables. *Michigan Math. J.* **30** 69–81. MR0694930

[7] CHOW, Y. S. and TEICHER, H. (1988). *Probability Theory*, 2nd ed. Springer, New York. MR0953964

[8] CSÖRGÖ, M. and HORVÁTH, L. (1997). *Limit Theorems in Change-Point Analysis.* Wiley, New York.

[9] DEDECKER, J. and DOUKHAN, P. (2003). A new covariance inequality and applications. *Stochastic Process. Appl.* **106** 63–80. MR1983043

[10] DEDECKER, J. and MERLEVÈDE, F. (2002). Necessary and sufficient conditions for the conditional central limit theorem. *Ann. Probab.* **30** 1044–1081. MR1920101

[11] DEDECKER, J. and MERLEVÈDE, F. (2003). The conditional central limit theorem in Hilbert spaces. *Stochastic Process. Appl.* **108** 229–262. MR2019054

[12] DEDECKER, J. and MERLEVÈDE, F. (2003). Convergence rates in the law of large numbers for Banach valued dependent variables. Technical Report LSTA 2003-4.

[13] DEDECKER, J. and PRIEUR, C. (2004). Coupling for $\tau$-dependent sequences and applications. *J. Theoret. Probab.* **17** 861–885. MR2105738

[14] DIACONIS, P. and FREEDMAN, D. (1999). Iterated random functions. *SIAM Rev.* **41** 45–76. MR1669737

[15] DOOB, J. (1953). *Stochastic Processes.* Wiley, New York. MR0058896

[16] DOUKHAN, P. (2003). Models, inequalities, and limit theorems for stationary sequences. In *Theory and Applications of Long-Range Dependence* (P. Doukhan, G. Oppenheim and M. S. Taqqu, eds.) 43–100. Birkhäuser, Boston. MR1956044

[17] EBERLEIN, E. (1986). On strong invariance principles under dependence assumptions. *Ann. Probab.* **14** 260–270. MR0815969

[18] EBERLEIN, E. and TAQQU, M. S., eds. (1986). *Dependence in Probability and Statistics: A Survey of Recent Results.* Birkhäuser, Boston. MR0899982

[19] ELTON, J. H. (1990). A multiplicative ergodic theorem for Lipschitz maps. *Stochastic Process. Appl.* **34** 39–47. MR1039561

[20] FELLER, W. (1971). *An Introduction to Probability Theory and Its Applications* **II**. Wiley, New York. MR0270403

[21] GORDIN, M. I. (1969). The central limit theorem for stationary processes. *Dokl. Akad. Nauk SSSR* **188** 739–741. MR0251785

[22] GORDIN, M. I. and LIFSIC, B. (1978). The central limit theorem for stationary Markov processes. *Doklady* **19** 392–394. MR0501277

[23] HALL, P. and HEYDE, C. C. (1980). *Martingale Limit Theory and Its Application*. Academic Press, New York.

[24] HANNAN, E. J. (1979). The central limit theorem for time series regression. *Stochastic Process. Appl.* **9** 281–289. MR0562049

[25] HAUSER, M. A. and KUNST, R. M. (2001). Forecasting high-frequency financial data with the ARFIMA–ARCH model. *J. Forecasting* **20** 501–518.

[26] HEYDE, C. C. (1975). On the central limit theorem and iterated logarithm law for stationary processes. *Bull. Austral. Math. Soc.* **12** 1–8. MR0372954

[27] HEYDE, C. C. and SCOTT, D. J. (1973). Invariance principles for the law of the iterated logarithm for martingales and processes with stationary increments. *Ann. Probab.* **1** 428–436. MR0353403

[28] HSING, T. and WU, W. B. (2004). On weighted $U$-statistics for stationary processes. *Ann. Probab.* **32** 1600–1631. MR2060311

[29] IBRAGIMOV, I. A. and LINNIK, YU. V. (1971). *Independent and Stationary Sequences of Random Variables*. Wolters-Noordhoff, Groningen. MR0322926

[30] KOMLÓS, J., MAJOR, P. and TUSNÁDY, G. (1975). An approximation of partial sums of independent RV's and the sample DF. I. *Z. Wahrsch. Verw. Gebiete* **32** 111–131. MR0375412

[31] KOMLÓS, J., MAJOR, P. and TUSNÁDY, G. (1976). An approximation of partial sums of independent RV's and the sample DF. II. *Z. Wahrsch. Verw. Gebiete* **34** 33–58. MR0402883

[32] KRENGEL, U. (1985). *Ergodic Theorems*. de Gruyter, Berlin. MR0797411

[33] LAI, T. L. and STOUT, W. (1980). Limit theorems for sums of dependent random variables. *Z. Wahrsch. Verw. Gebiete* **51** 1–14. MR0566103

[34] LIEN, D. and TSE, Y. K. (1999). Forecasting the Nikkei spot index with fractional cointegration. *J. Forecasting* **18** 259–273.

[35] LIN, Z. and LU, C. (1996). *Limit Theory for Mixing Dependent Random Variables*. Kluwer, Dordrecht. MR1486580

[36] MCLEISH, D. L. (1975). A maximal inequality and dependent strong laws. *Ann. Probab.* **3** 829–839. MR0400382

[37] MENCHOFF, D. (1923). Sur les series de fonctions orthogonales I. *Fund. Math.* **4** 82–105.

[38] MERLEVÈDE, F. and PELIGRAD, M. (2007). On the central limit theorem and its weak invariance principle under projective criteria. *J. Theoret. Probab.* To appear.

[39] MORICZ, F. (1976). Moment inequalities and strong laws of large numbers. *Z. Wahrsch. Verw. Gebiete* **35** 299–314. MR0407950

[40] PELIGRAD, M. and UTEV, S. (2005). A new maximal inequality and invariance principle for stationary sequences. *Ann. Probab.* **33** 798–815. MR2123210

[41] PELIGRAD, M., UTEV, S. and WU, W. B. (2007). A maximal $\mathbb{L}_p$-inequality for stationary sequences and its applications. *Proc. Amer. Math. Soc.* **135** 541–550. MR2255301

[42] PHILIPP, W. (1969). The law of the iterated logarithm for mixing stochastic processes. *Ann. Math. Statist.* **40** 1985–1991. MR0317390

[43] Philipp, W. (1986). Invariance principles for independent and weakly dependent random variables. In *Dependence in Probability and Statistics: A Survey of Recent Results* (E. Eberlein and M. S. Taqqu, eds.) 225–268. Birkhäuser, Boston. MR0899992

[44] Philipp, W. and Stout, W. (1975). Almost sure invariance principles for partial sums of weakly dependent random variables. *Mem. Amer. Math. Soc.* **2**. MR0433597

[45] Phillips, P. C. B. (1999). Unit root log periodogram regression. Technical Report 1244, Cowles Foundation, Yale Univ.

[46] Priestley, M. B. (1988). *Nonlinear and Nonstationary Time Series Analysis.* Academic Press, London. MR0991969

[47] Reznik, M. H. (1968). The law of the iterated logarithm for certain classes of stationary processes. *Theory Probab. Appl.* **13** 606–621. MR0243592

[48] Rio, E. (1995). The functional law of the iterated logarithm for stationary strongly mixing sequences. *Ann. Probab.* **23** 1188–1203. MR1349167

[49] Romano, J. P. and Thombs, L. A. (1996). Inference for autocorrelations under weak assumptions. *J. Amer. Statist. Assoc.* **91** 590–600. MR1395728

[50] Serfling, R. J. (1970). Moment inequalities for the maximum cumulative sum. *Ann. Math. Statist.* **41** 1227–1234. MR0268938

[51] Shao, Q. M. (1993). Almost sure invariance principles for mixing sequences of random variables. *Stochastic Process. Appl.* **48** 319–334. MR1244549

[52] Shao, X. and Wu, W. B. (2007). Local Whittle estimation of fractional integration for nonlinear processes. *Econometric Theory.* To appear.

[53] Stine, R. A. (2006). Nonlinear time series. In *Encyclopedia of Statistical Sciences*, 2nd ed. (S. Kotz, C. B. Read, N. Balakrishnan and B. Vidakovic, eds.) 5581–5588. Wiley, New York.

[54] Stout, W. F. (1970). The Hartman–Wintner law of the iterated logarithm for martingales. *Ann. Math. Statist.* **41** 2158–2160.

[55] Stout, W. F. (1974). *Almost Sure Convergence.* Academic Press, New York. MR0455094

[56] Strassen, V. (1964). An invariance principle for the law of the iterated logarithm. *Z. Wahrsch. Verw. Gebiete* **3** 211–226. MR0175194

[57] Strassen, V. (1967). Almost sure behaviour of sums of independent random variables and martingales. *Proceedings of the Fifth Berkeley Symposium of Math. Statist. Probab.* **2** 315–343. Univ. California Press, Berkeley. MR0214118

[58] Tien, N. D. and Huang N. V. (1989). On the convergence of weighted sums of martingale differences. *Probability Theory on Vector Spaces IV. Lecture Notes in Math.* **1391** 293–307. Springer, Berlin. MR1020569

[59] Tong, H. (1990). *Nonlinear Time Series. A Dynamical System Approach.* Oxford Univ. Press. MR1079320

[60] Volný, D. (1993). Approximating martingales and the central limit theorem for strictly stationary processes. *Stochastic Process. Appl.* **44** 41–74. MR1198662

[61] Volný, D. and Samek, P. (2000). On the invariance principle and the law of the iterated logarithm for stationary processes. In *Mathematical Physics and Stochastic Analysis (Lisbon, 1998)* 424–438. World Scientific, River Edge, NJ. MR1893125

[62] Woodroofe, M. (1992). A central limit theorem for functions of a Markov chain with applications to shifts. *Stochastic Process. Appl.* **41** 33–44. MR1162717

[63] Woyczyński, W. A. (1982). Asymptotic behavior of martingales in Banach spaces. II. *Martingale Theory in Harmonic Analysis and Banach Spaces. Lecture Notes in Math.* **939** 216–225. Springer, Berlin–New York. MR0668549

[64] Wu, W. B. (2005). Nonlinear system theory: Another look at dependence. *Proc. Natl. Acad. Sci. USA* **102** 14150–14154. MR2172215

[65] Wu, W. B. (2005). On the Bahadur representation of sample quantiles for dependent sequences. *Ann. Statist.* **33** 1934–1963. MR2166566

[66] Wu, W. B. and Min, W. (2005). On linear processes with dependent innovations. *Stochastic Process. Appl.* **115** 939–958. MR2138809

[67] Wu, W. B. and Shao, X. (2004). Limit theorems for iterated random functions. *J. Appl. Probab.* **41** 425–436. MR2052582

[68] Wu, W. B. and Woodroofe, M. (2004). Martingale approximations for sums of stationary processes. *Ann. Probab.* **32** 1674–1690. MR2060314

[69] Wu, W. B. and Zhao, Z. (2007). Inference of trends in time series. *J. Roy. Statist. Soc. Ser. B.* To appear. Available at http://galton.uchicago.edu/˜wbwu/papers/bard-nov1-06.pdf.

[70] Wu, W. B. and Zhao, Z. (2006). Moderate deviations for stationary processes. *Statist. Sinica.* To appear. Available at http://galton.uchicago.edu/techreports/tr571.pdf.

[71] Yokoyama, R. (1995). On the central limit theorem and law of the iterated logarithm for stationary processes with applications to linear processes. *Stochastic Process. Appl.* **59** 343–351. MR1357660

[72] Zhao, O. and Woodroofe, M. (2006). Law of the iterated logarithm for stationary processes. Technical report. Available at http://www.stat.lsa.umich.edu/˜michaelw/PPRS/lilsbmt.pdf.

Department of Statistics
University of Chicago
5734 S. University Avenue
Chicago, Illinois 60637
USA
E-mail: wbwu@galton.uchicago.edu